\documentclass[12pt]{article}
\usepackage{philstyle}

\usepackage{xypic,amsmath,graphicx,epsfig,psfrag}

\def\cF{{\cal F}}

\def\QED{\ $\blacksquare$\smallskip}

\def\de#1{\textit{#1}}
\def\circle{{S^1}}

\def\complexes{{\mathbb C}}
\def\rationals{{\mathbb Q}}
\def\reals{{\mathbb R}}
\def\integers{{\mathbb Z}}

\def\raw{\rightarrow}

\def\I{^{-1}}

\def\tD{D^{(\infty)}}

\def\tf{\tilde{f}}
\def\th{\tilde{h}}

\def\RBSk{S^{(k)}}
\def\RBSkC{S^{(k)}(\C)}

\def\Rk{R^{(k)}}

\def\RBSkC{S^{(k)}_{\C}}
\def\RkC{R^{(k)}_{\C}}

\def\htop{h_{top}}
\def\homeo{homeomorphism}
\def\homeos{homeomorphisms}

\def\st{|}

\def\isomorphic{\cong}
\def\Z{\integers}
\def\C{\complexes}

\def\R{\reals}

\DeclareMathOperator{\id}{id}
\DeclareMathOperator{\im}{im}

\def\tGamma{\tilde{\Gamma}}
\def\tx{\tilde{x}}

\def\pA{pseudo-Anosov}

\def\op{\orient{p}}
\def\opstar{\orient{p_*}}

\def\ophistar{\orient{\phi_*}}

\def\tPhit{{\Phi}^{(2)}}
\def\tGamma{\tilde{\Gamma}}
\def\tX{\tilde{X}}

\renewcommand{\exp}[1]{e^{#1}}

\newcommand{\bd}{\partial}
\newcommand{\orient}[1]{\overrightarrow{#1}}

\newcommand{\tk}[1]{{#1}^{(k)}}
\newcommand{\tinf}[1]{{#1}^{(\infty)}}

\DeclareMathOperator{\lcm}{lcm}
\DeclareMathOperator{\spec}{sr}
\DeclareMathOperator{\GL}{GL}
\DeclareMathOperator{\Mat}{Mat}

\newcommand{\red}{\Phi}
\DeclareMathOperator{\rank}{rank}

\begin{document}
\title{The Burau matrix and the entropy of a braid}

\title{Estimating a braid's entropy using the Burau matrix}

\title{The Burau estimate for the entropy of a braid}

\author{Gavin Band and Philip Boyland}

{\namefont

\maketitle

}




\begin{abstract}
  The topological entropy of a braid is the infimum of the entropies
  of all homeomorphisms of the disk which have a finite invariant set
  represented by the braid. When the isotopy class represented by the
  braid is pseudo-Anosov or is reducible with a pseudo-Anosov component, this
  entropy is positive.  Fried and Kolev proved that the entropy is
  bounded below by the logarithm of the spectral radius of the braid's
  Burau matrix, $B(t)$, after substituting a complex number of
  modulus~$1$ in place of $t$. In this paper we show that for a
  pseudo-Anosov braid the estimate is sharp for the substitution of a
  root of unity if and only if it is sharp for $t=-1$. Further, this
  happens if and only if the invariant foliations of the pseudo-Anosov
  map have odd order singularities at the strings of the braid and all
  interior singularities have even order. An analogous theorem for
  reducible braids is also proved.
\end{abstract}

\section{Introduction}
Artin's braid group and its Burau representation have been extensively
studied by many researchers from many points of view. In dynamical
applications a braid is often used to describe the motion of a
collection of points in the two-dimensional disk. Since the braid
depends only on the motion of the points, it is describing an isotopy
class of \homeos\ on the complement of the points. Thus, the
interpretation of the braid group on $n$-strings, $B_n$, as a mapping
class group of the $n$-punctured disk is frequently used, and so
Thurston's classification theorem for surface isotopy classes is an
important tool.

The (reduced) Burau matrix, $B(t)$, of a braid $\beta\in B_n$, is an
$(n-1) \times (n-1)$ matrix with entries in $\integers[t, t\I]$, i.e.\
the entries of the matrix are Laurent polynomials over the integers.
In the early 1980's two different but closely related dynamical
interpretations of the Burau matrix emerged.  Using the construction
in Franks' paper \cite{Franks1981}, the Burau matrix can be
interpreted as the signed, linking matrix of a certain Axiom A flow
associated with the braid. The signed, linking matrix is an enhanced
Markov transition matrix which records the linking of the Markov boxes
with the strings of the braid as well as the orientations of their
images.

The second dynamical interpretation comes from the machinery in
Fried's paper \cite{Fried1986}.  In this case the Burau matrix of a
braid arises as the induced action of the associated mapping class on
a $\integers$-cover, where the first homology of the cover is given
the structure of a module over $\integers[t, t\I]$ (this is a
standard description, see, for example, \cite{BirmanAndBrendle2005}).
In the more general setting of twisted cohomology, Fried showed that
after using the appropriate representation of the fundamental group,
the spectral radius of the induced action gives a lower bound on the
topological entropy. While these interpretations of the Burau matrix
were not explicit in either of these two papers, the two authors were
certainly aware of them (P.B. personal communication, 1984).  See
\cite{BoylandWarwick} for an introductory exposition of these two
interpretations.

Because the topological entropy of a self-map measures a certain kind
of exponential growth it is natural to expect that, at least in
certain cases, it is detectable from the growth rates of induced maps
on various algebraic objects associated with the underlying
space. Thus, for example, the growth rate of induced map on first
homology (i.e.\ its spectral radius) gives a lower bound for the
topological entropy (\cite{Manning1974}), as does exponential growth
rate of word length in the fundamental group under iteration by $f_*$
(\cite{Bowen78, Orsay1979}).

These lower bounds only depend on the homotopy class of the map, and
so it is also natural to ask whether the bounds are attained, i.e.\ is
there a map in the homotopy class that realises the lower bound? For
surface \homeos\ this question was answered by Thurston. One
consequence of his classification theorem is that any isotopy class of
surface \homeos\ contains a map with entropy equal to the growth rate
on the fundamental group (\cite{Thurston1988, Orsay1979}).  While this
result is invaluable in theory, in practice, the computation of word
length growth in non-Abelian groups is very difficult.  On the other
hand, computations in homology are much more tractable, but frequently
give only trivial lower bounds.  A fundamental idea in Fried's paper
\cite{Fried1986} is that there is a middle ground between these two
theories provided by the action on twisted homology. In the most
concrete case this amounts to providing a systematic way to examine
the growth rate of the action on homology in a collection of finite
covers.

Thus one sees that the Burau representation provides a lower bound for
the topological entropy of the isotopy class represented by a
braid. Specifically, if $h$ is a \homeo\ of the $n$-punctured disk
which is represented by the braid $\beta\in B_n$ with Burau matrix
$B(t)$, then
\begin{equation}\label{firstburest}
  \htop(h) \geq \sup\{\log\spec(B(\eta)): \eta\in\circle\}
\end{equation}
where $\spec(B(\eta))$ is the spectral radius of the complex
matrix~$B(\eta)$ obtained by substituting the complex number $\eta$
with $|\eta| = 1$ into the Burau matrix.  This estimate was obtained
directly with different methods by Kolev~\cite{Kolev1989}. The
estimate in \eqref{firstburest} and its alternative version in
\eqref{burest} below will be called the \de{Burau estimate}.  If the
inequality in \eqref{firstburest} is an equality for $\eta = \eta_0$,
then the Burau estimate is said to be \de{sharp at~$\eta_0$}.

Since a braid represents an isotopy class, we may use Thurston's
classification scheme to classify braids. Thus a braid is said to be
\pA, finite order or reducible if its corresponding isotopy class
is. If a braid is finite order or reducible with all finite order
components, there is a map in the class with zero topological entropy,
and so the Burau estimate is already sharp at~$\eta = 1$.  The main
result here for \pA\ braids is
\begin{theorem}\label{introthm1}
  For a \pA\ braid $\beta$, the Burau estimate is sharp at some root
  of unity~$\eta_0$ if and only if it is sharp at~$-1$. Further,
  this happens if and only if the invariant foliations for a \pA-map
  in the class represented by $\beta$ have odd order singularities at
  all punctures and all interior singularities are even order.
\end{theorem}

A portion of this theorem was obtained in \cite{Song2002} (see
Remark~\ref{rk:song} below).  The substitution of complex numbers on
the unit circle which are not roots of unity requires different
methods. In a subsequent paper we will show that for a \pA-braid,
$\spec(B(e^{2 \pi i \theta})) < \lambda$ for all $\theta\not\in
\rationals$.

There is an analogous theorem for reducible braids with at least one
\pA\ component. Its statement is rather complicated (see
Theorem~\ref{theorem:spectral radius for reducible braids} below), but
one useful consequence is the following.

\begin{theorem}\label{introthm2b}
  For any braid $\beta$ on $n$ strings, the Burau estimate is sharp at
  some root of unity~$\eta_0$ if and only it is sharp at a~$k$-th root
  of unity~$e^{2 \pi i j/k}$, for some~$k \leq \tfrac{2}{3} n$ which
  is a power of~$2$.
\end{theorem}

There are two main components in the proof of these theorems.  The
main algebraic tool is Lemma~\ref{eiglem} which shows that the
union of the spectra obtained by substituting all the $k^{th}$ roots
of unity into the Burau matrix $B(t)$ yields essentially the entire
spectrum of the action on first homology of a lift of the
corresponding mapping class to the $k$-fold cover.  This is coupled
with information about the connection of entropy, the Thurston normal
form, the action on homology, and the orientability of the invariant
foliations of a \pA-map.

The investigations of this paper were partly inspired by recent work using the
Burau estimate to get lower bounds on the entropy of fluid flows by applying the
Burau estimate to the braids generated by large collections of points moving
with the fluid (\cite{Thiffeault2005,Gouillart2006,Thiffeault2006}). 
For these applications a good understanding of
``sharpness'' of the estimate is necessary and this paper provides a first
step. In addition, questions surrounding the Burau representation provide an
important, special case of the more general question of dynamics on abelian
covers of surfaces which we investigate in a subsequent paper.

\section{Preliminaries}\label{sect:prelim}
\subsection{Standing hypotheses and conventions}
\label{sect:standing}
In this paper surfaces $X$ will always be orientable, 
perhaps with boundary, and compact except
perhaps for a finite number of punctures. We fix a Riemannian metric on $X$
which allows us to speak of the lengths of tangent vectors. Self-homeomorphisms
of the surface $f:X\raw X$ are always orientation-preserving. 
If no coefficient ring for homology is specified, it is 
assumed  to be the integers $\Z$,  and so $H_1(X) := H_1(X;\Z)$.
The induced map of the \homeo\ $f$ on first homology is denoted
$f_*$.  If $M$ is a square, complex matrix, then $\spec(M)$ denotes 
its spectral radius.

The classification theorem for regular connected covering spaces
identifies each such cover with a normal subgroup of the fundamental
group of the connected base space~$X$, or equivalently, with an epimorphism
$\pi_1(X)\raw G$, where $G$ is the group of deck transformations of
the cover (see, for example, \cite{FultonAlgebraicTopology}).   In this paper $G$ will always be
abelian, and so we often designate a cover $\tX$ by an epimorphism
$\rho:H_1(X) \raw G$, with the Hurewicz homomorphism $\pi_1(X) \raw
H_1(X)$ being implicit.

More generally, we shall also need to consider disconnected covers
over connected and disconnected base spaces. In these cases it will be
convenient to continue to designate the cover by a homomorphism
$\rho:H_1(X) \raw G$, which perhaps is not surjective.  As this is
less commonly encountered, we describe the cover~$\tilde{X}$
associated to such a homomorphism~$\rho$.  First suppose~$X$ is
connected, and let~$G' = \im \rho \subset G$.  Then~$\rho$ induces an
epimorphism~$\rho':H_1(X) \raw G'$ with the same domain as~$\rho$ but
with range~$G'$, and so as above it determines a connected covering
space~$\tilde{X}_0$ over~$X$ with deck group~$G'$.  We
define~$\tilde{X}$ to be the disjoint union of copies
of~$\tilde{X}_0$, one such copy for each coset of~$G'$ in~$G$.  The
action of~$G'$ by deck transformations on each copy of~$\tilde{X}_0$
extends to an action of~$G$ on~$\tilde{X}$, in such a way that for
any~$g, g' \in G$, the element~$g$ sends the copy of~$\tilde{X}_0$
corresponding to~$g'+G$ to that corresponding to~$g+g'+G'$.
Then~$\tilde{X} / G \isomorphic X$, and so this makes~$\tilde{X}$ a
covering space over~$X$ with deck group~$G$.  Finally, if~$X$ is
disconnected, we define~$\tilde{X}$ to be the disjoint union, over all
connected components~$Y$ of~$X$, of the covering space~$\tilde{Y}$
corresponding to~$\rho|_{H_1(Y)}$ as above.

In all cases every \homeo\ of the base $f:X\raw X$ will satisfy $\rho
f_* = \rho$, which implies that $f$ lifts to $\tf:\tX\raw\tX$ which
commutes with all deck transformations $g\in G$.

\subsection{The Braid group and the Burau Representation}\label{sub:braid}
In this section we briefly survey relevant results from the point of
view of Dynamical Systems.  The classic reference
are~\cite{Artin1925,Birman1974}, and see \cite{BirmanAndBrendle2005}
for a survey including recent developments and \cite{Boyland1994} for
a survey of dynamical applications.  The book \cite{Fenn1983} and the
classic paper \cite{Milnor1968} are good sources of information on
homology of cyclic covers.

The braid group on $n$~strings, $B_n$, is defined 
using generators and relations as
\begin{equation*}
  B_n = \left\langle \sigma_1, \dots, \sigma_{n-1} |
    \text{$\sigma_i \sigma_k = \sigma_k \sigma_i$ if~$|i-k|>2$,
      $\sigma_i\sigma_{i+1}\sigma_i=\sigma_{i+1}\sigma_i\sigma_{i+1}$} 
  \right\rangle.
\end{equation*}
In this paper we shall be principally concerned with $B_n$ interpreted
as a mapping class group, namely, the group of isotopy classes of
homeomorphisms of the~$n$-punctured disk where all homeomorphisms and
all isotopies are required to fix the boundary pointwise.  Letting
$x_j = j/(n+1)$ for $j = 1, \dots, n$ and $D_n=\{z \in\C : |z-1/2|\leq
1/2\} \setminus \{x_1, \dots, x_n\}$, the generator~$\sigma_i$ of
$B_n$ corresponds to a \homeo\ that switches $x_i$ and $x_{i+1}$ in a
counter-clockwise direction.  When we indicate a braid $\beta\in B_n$
we will always be identifying it with its corresponding isotopy class,
and so for example, $h\in\beta\in B_n$ means that $h$ is a \homeo\ of $D_n$
that is contained in the isotopy class corresponding to $\beta$.

As with the braid group we shall need an interpretation of the Burau
representation with dynamical content, as the action on homology in a
particular cover.
To construct the $\integers$-cover relevant to the Burau
representation, we fix a basepoint~$x_0 \in D_n$, and around each
puncture~$x_i$ we take a small clockwise loop~$\Gamma_i$ which we then
homotope so it begins and ends at~$x_0$.  Then $\pi_1(D_n, x_0)$ and
$H_1(D_n)$ are freely generated by the set~$\{\gamma_i\}$ of
homotopy/homology classes of the~$\Gamma_i$'s.  Let~$\tau$ be the
epimorphism of~$H_1(D_n)$ onto~$\Z$ generated by~$\tau(\gamma_i) = 1$
for all $i$; the resulting cover is called the \de{Burau cover} and is
denoted $\tD$. By construction the deck group of $\tD$ is isomorphic
to $\Z$, and we call its generator $T$.  An orientation-preserving
\homeo~$h$ of~$D_n$ must permute the punctures of $D_n$ and will
therefore act as a permutation on the generators of~$H_1(D_n)$, and so
$\tau h_* = \tau$.  Thus any orientation-preserving \homeo~$h$
of~$D_n$ will lift to a \homeo~$\th:\tD\raw\tD$ and further, each lift
of $h$ commutes with all deck transformations, or $T^j \th = \th T^j$
for all $j\in\integers$.  Note that by definition, any~$h \in \beta\in
B_n$ is a \homeo\ that fixes the outside boundary of~$D_n$ point-wise
and this yields a preferred lift of $h$ to $\tD$, namely, the lift
which fixes the lift of the outside boundary of $D_n$
point-wise. Unless indicated otherwise any lift to $\tD$ of an
$h\in\beta$ will be this preferred lift.

Next note that the first homology group $H_1(\tD)$ has a natural
structure as a module over~$R := \integers[t^{\pm 1}]$, the ring of
all Laurent polynomials with coefficients in $\Z$, i.e.\ the group ring
of $\Z$.  To describe this structure, lift all the loops $\Gamma_i$ to
arcs $\tGamma_i\in\tD$ with all $\tGamma_i$ starting at some point
$\tx_0$ and ending at the point $T \tx_0$. Thus for $i = 1, \dots,
n-1$, $\zeta_i := [\tilde{\Gamma}_{i+1} - \tilde{\Gamma}_{i}]\in
H_1(\tD)$.  For a Laurent polynomial $p(t) = \sum a_j t^j \in R$,
$p(t)\zeta_i$ represents the class $\sum a_j T^j \zeta_i \in
H_1(\tD)$, and so as an $R$-module, $H_1(\tD)\isomorphic R^{n-1}$.
Now since the lift of a \homeo, $\th$, commutes with~$T$, we see that
that~$\th$ acts on~$H_1(\tD)$ by an~$R$-module isomorphism.  So with
respect to the $R$-module basis~$\{\zeta_i\}$ of~$H_1(\tD)$, $\th$
acts by a matrix~$B(\th) \in \GL(n-1,R)$.

For a braid~$\beta \in B_n$, pick~$h \in \beta$ and its preferred
lift~$\th$ to~$\tD$.  The matrix~$B(\beta) = B(\th)$ is called the
\emph{reduced Burau matrix} of~$\beta$, and the corresponding
homomorphism\linebreak~$B_n\raw\GL(n-1,R)$ the \emph{reduced Burau
 representation} of the braid group~$B_n$. When the braid $\beta$
is fixed, we often will write its Burau matrix as $B(t)$.
The full Burau
representation will not be used here, but for completeness we note
that it can be defined similarly using the action of~$\th$ on the
relative homology group~$H_1(\tD,F)$, where~$F$ is the fiber above the
basepoint~$x_0$ (which in this case must be a fixed point of~$h$).

Since our results work for all $n>2$, we fix once and for all such an
$n$ and suppress the dependence of objects on $n$ when possible.  So,
for example, we write just $D$ not $D_n$.

\subsection{The Nielsen-Thurston normal form}\label{sub:niel}
Since we have identified the braid group $B_n$ with a surface mapping
class group, Thurston's classification theorem will be of central
importance here. This theorem identifies ``simplest'' representatives
in any isotopy class.  See \cite{Orsay1979} and \cite{Thurston1988}
for more information. There are minor differences in the literature in
how punctures, boundary and the reducible case are handled in stating
Thurston's results. The version we give in Theorem~\ref{thurstons
  theorem} is adapted to our use with the braid group.

The two main ingredients in Thurston's classification are finite order
and \pA\ maps. A map $\phi$ is \de{finite order} if $\phi^n = \id$ for
some $n \geq 1$.  The map $\phi$ is \pA\ if there are a pair of
transverse, $\phi$-invariant measured foliations\footnote{The
  designation ``measured foliation'' is the standard shortening of the
  more proper, and much longer name, ``foliation with conical
  singularities with a holonomy invariant transverse measure''.},
$\cF^u$ and $\cF^s$. Under the action of $\phi$ the transverse
measures are expanded and contracted by a number $\lambda > 1$, which
is called the \de{expansion factor} of the \pA\ map. This fact is
usually indicated by the notation $\phi_* \cF^u = \lambda \cF^u$ and
$\phi\I_* \cF^s = \lambda \cF^s$.  The supporting surface of a
pseudo-Anosov map~$\phi$ may be connected or disconnected, but in the
latter case we require that~$\phi$ cyclically permutes the components.

\begin{figure}
  \centering\epsfig{figure=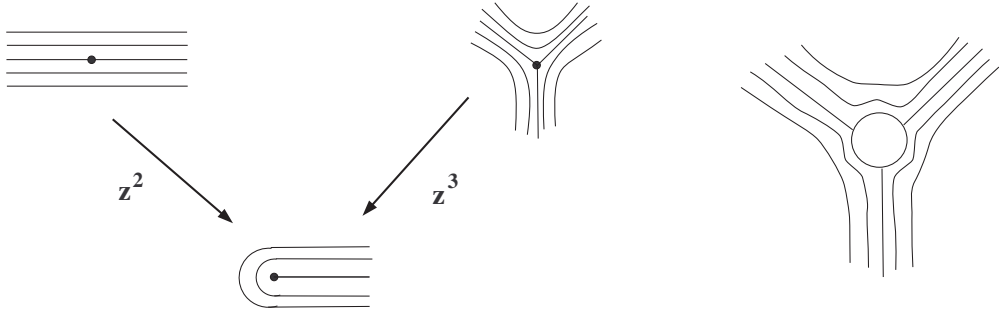,width=.80\linewidth}
  \caption{(a) Constructing prongs; (b) a boundary three prong.}
  \label{singfig}
\end{figure}

Since the structure of the singular points in a measured foliation is
central here we describe it in more detail.  Near a \de{regular point}
a measured foliation looks locally like the foliation of $\reals^2$ by
horizontal lines.  A measured foliation is also allowed to have a
finite number of non-regular or singular points which are required to
have a very a specific local structure which is characterized by the
order of the singularity, i.e.\  by the number of leaves coming
directly into the point. The local structure can be succinctly
described using covers branched over ${\bf 0} \in\complexes$. Starting
with the foliation of $\C$ by horizontal lines and projecting under
$z\mapsto z^2$, we get the local structure of an \de{order one or
  one-prong singularity} at ${\bf 0}$ (see Figure~\ref{singfig}(a)).
  Lifting the one-prong by the map $z \mapsto z^n$ for $n>2$
gives an \de{order $n$ or $n$-prong singularity}. An \de{order
  $n$-punctured singularity} is formed by removing the singular point
from an $n$-prong, and an \de{order $n$-boundary singularity} is
obtained by replacing this puncture with a boundary circle (see
Figure~\ref{singfig}(b)).

The invariant foliations associated with a \pA\ map have a few special
additional qualifications.  They can have punctured or boundary
one-prongs, but interior one-prongs do not occur as they will not
persist under isotopy. An interior two-prong is a regular point and is
not considered a singularity, but a punctured or boundary two-prong is
considered a singularity. Given a measured foliation $\cF$ on a
surface $X$ with genus $g$, the Euler-Poincar\'e-Hopf formula says
that
\begin{equation}
  \label{EPHFormula}
  2 -2g = \sum (1 - \kappa_i/2),
\end{equation}
where the sum is over singularities of all types (interior,
boundary and punctured) of $\cF$, and $\kappa_i$ is the order of the
$i^{th}$ singularity. 

Now recall that any \homeo\ $h\in\beta\in B_n$ must
fix the boundary of $D_n$ pointwise, a property
that  is not shared by the Thurston representatives under
the usual definitions. 
In addition, the isotopies used in the braid group must fix the
boundary pointwise. These two facts require a small alteration
in the designation of Thurston representative in an isotopy class. 
Here is a version 
of Thurston's classification theorem adapted to our situation.
It follows, for example, by altering the constructions in
\cite{Boyland1999}.

\begin{theorem}[Thurston] 
  \label{thurstons theorem}
  Let~$f$ be a homeomorphism of the possibly disconnected, possibly
  bordered surface~$X$, compact except for a finite number of
  punctures, and assume that~$f$ fixes the boundary of~$X$ pointwise.
  Then there is a homeomorphism~$\Phi$ isotopic to~$f$ by an isotopy
  which fixes the boundary pointwise, and a decomposition
  \begin{equation*}
    X = A \cup \bigcup_{j=1}^m X_j
  \end{equation*}
  of~$X$ into pairwise disjoint~$\Phi$-invariant sets, 
  with the following properties.
  \begin{enumerate}
  \item $A$ is a finite disjoint union of embedded open annuli. 
    The waist curve of such an annulus is never
    null-homotopic; nor are the waist curves of two annuli
    mutually homotopic.  If the waist curve of an
    annulus~$a \in A$ is homotopic to a boundary 
    component~$b$ of~$X$, then~$b$ is a component of~$\partial a$
    (and by convention we include~$b$ in~$a$).
  \item Each~$X_j$ is the union of a collection of connected
    components of~$X \setminus A$ which are permuted cyclically
    by~$\Phi$; and the restriction of~$\Phi$ to~$X_j$ is either finite
    order or \pA.
  \item The restriction of~$\Phi$ to~$A$ has zero topological entropy.
  \end{enumerate}
\end{theorem}

\subsection{Branched Covering spaces and oriented foliations}\label{sub:branch}
Given a finite set $B \subset X$, the triple $p:Y\raw X$ is a called
\de{a covering space branched over} $B$ if $p$ is onto and restricts
to a covering map (i.e.\ a surjective local homeomorphism which is
evenly covered over any small neighborhood,
c.f. \cite{FultonAlgebraicTopology}) of $Y \setminus p^{-1} B$ over $X
\setminus B$.  Now let~$\cF$ be a measured foliation on the
surface~$X$.  The foliation~$\cF$ is \emph{orientable} provided there
is a vector field~$\Upsilon$ on~$X$, zero at the singularities
of~$\cF$, and everywhere else nonzero and tangent to~$\cF$.
When~$\cF$ is not orientable, it is often useful to orient it by
lifting it to a two-sheeted branched cover constructed as follows.
Start by puncturing $X$ at the singularities of the foliation and
denote the resulting space by $X'$.  Now define $\orient{X'}$ as the
space of unit tangent vectors to~$\cF$ in $X'$ with the topology
induced as a subspace of the unit tangent bundle of $X$.  It is
evident that $\orient{X'}$ is a two-sheeted cover over~$X'$.  By
sewing back in the punctures, we obtain a two-fold branched cover
$\orient{p}:\orient{X}\raw X$, and pulling back $\cF$, we obtain a
measured foliation~$\orient{\cF}$ on~$\orient{X}$.  The
space~$\orient{X}$ equipped with the pulled back
foliation~$\orient{\cF}$ is called the \emph{orientation double-cover}
of the non-orientable foliation~$\cF$.

If $\gamma$ is a smooth loop in $X'$, then since $X$ is an orientable
surface, the tangent bundle restricted to $\gamma$ is trivial and so
there is a well defined monodromy as we pull along $\gamma$ the unit
vectors tangent to $\cF$. If this monodromy brings a vector back to
itself, we say that \de{the foliation is oriented along $\gamma$}, and
if it brings a vector back to its opposite, \de{the foliation is
  disoriented along $\gamma$}. It is evident from the definition of
the orientation double cover that in the first case $\gamma$ lifts to
a pair of disjoint loops in $\orient{X}$, while in the latter case,
$\gamma$ is covered by a single loop $\tilde{\gamma}\subset X'$, and
the covering map induces a degree-two map
$\tilde{\gamma}\raw\gamma$. Thus a foliation is oriented if and only
if it is oriented along every loop in the complement of the
singularity set, and if $\orient{\rho}:H_1(X')\raw \Z_2$ is the
epimorphism associated to the cover $\orient{X'}$, then the foliation
is oriented along $\gamma$ if and only if $\orient{\rho}([\gamma]) =
0$.

The next lemma gives a simple criterion for checking when the
foliation is oriented when lifted to a cover.  Its proof is standard
and we omit it. Implicit in the statement of (ii) is the fact that if
$\gamma$ is a small loop surrounding a singularity $P$, then
$\orient{\rho}([\gamma]) = 0$, if $P$ has even order, and
$\orient{\rho}([\gamma]) = 1$, if $P$ has odd order.  Thus if all
interior singularities of a measured foliation on $X$ are of even
order, then we may treat $\orient{\rho}$ as being defined on $H_1(X)$.

\begin{lemma}\label{orientlemma}
  Assume that~$X$ is a possibly disconnected, possibly bordered
  surface~$X$, compact except for finitely many punctures, and
  that~$\cF$ is a measured foliation on $X$.  Let~$\rho:H_1(X)\raw G$
  be a homomorphism of~$H_1(X)$ to a finite abelian group~$G$ and
  let~$\tilde{X}$ be the corresponding covering space of~$X$.  The
  following are equivalent:
  \begin{enumerate}
  \item[(i)]~The lift of~$\cF$ to~$\tilde{X}$ is orientable;
  \item[(ii)]~All singularities in the interior of~$X$ have even order, and
    there exists a homomorphism~$\delta:\im\rho\raw\Z_2$ such
    that~$\delta\circ\rho=\orient{\rho}$, where $\orient{\rho}$ is
    the morphism defining the orientation cover of $\cF$.
  \end{enumerate}
\end{lemma}

\section{Finite Covers and substituting roots of unity}\label{sect:finite}

In this section we show that the substitutions of $k^{th}$ roots of
unity into the Burau representation of a braid give all the essential
spectral information about the action of the braid's mapping class on
homology in the $k$-fold cover.  The reducible case considered in
Section~\ref{sect:reducible} requires us to work with more general
subsurfaces of~$D$ and, in fact, the results of this section apply to
fairly general topological spaces.

Let~$h:X\raw X$ be a homeomorphism of the perhaps disconnected
space~$X$.  Suppose that~$\rho:H_1(X) \raw \Z$ is some homomorphism
which satisfies~$\rho h_* = \rho$, and denote
by\linebreak~$\tinf{p}:\tinf{X}\raw X$ the covering space associated to~$\rho$
as in Section \ref{sect:standing}.  The covering~$\tinf{p}$ is thus
generated by a deck group isomorphic to~$\Z$, and we denote the
generator of this deck group by~$T$.  Let~$\tinf{h}$ be a lift of~$h$
to~$\tinf{X}$.  The condition~$\rho h_* = \rho$ implies
that~$\tinf{h}$ commutes with~$T$ and hence also with every other
element of the deck group.

The main example for our purposes is where~$X = D$ is
the~$n$-punctured disk,~$h:D\raw D$ is a representative of the
braid~$\beta \in B_n$, and~$\rho = \tau$ is the homomorphism defining
the Burau cover~$\tinf{D}$ of~$D$ (see Section \ref{sub:braid}).  To
deal with reducible braids we will also have to consider
certain~$h$-invariant subsurfaces of~$D$ (which may be disconnected).
In general, all we shall assume about the space~$X$ is that its first
homology group~$H_1(X)$ is free and of finite rank.  The results can
easily be extended to work whenever the torsion-free part of~$H_1(X)$
has finite rank.

To construct cyclic covers, for each integer~$k > 0$,
let~$\xi_k:\Z\raw\Z_k$ be the quotient homomorphism and define~$\rho_k
= \xi_k \circ \rho$.  Let~$\tk{p}:\tk{X}\raw X$ be the covering space
associated to~$\rho_k$, so~$\tk{X} = \tinf{X} / T^k$.
Let~$\tk{q}:\tinf{X}\raw\tk{X}$ be the covering projection.  The image
under~$\tk{q}_*$ of~$H_1(\tinf{X})$ is a subgroup of~$H_1(\tk{X})$
which we denote by~$\RBSk$.

The map~$\tinf{h}$ pushes down to a well-defined lift~$\tk{h}$ of~$h$
on~$\tk{X}$.  Also, the deck group generator~$T$ of~$\tinf{X}$ pushes
down to a generator of the deck group for~$\tk{p}$, which we also
denote by~$T$.  Both~$\tk{h}_*$ and~$T_*$ leave~$\RBSk$ invariant by
definition.  

As in Section~\ref{sub:braid}, we write~$R = \Z[t^{\pm 1}]$ for the
ring of all Laurent polynomials with coefficients in~$\Z$ (i.e.\ the
group ring of~$\Z$) and we treat the integral homology
group~$H_1(\tinf{X})$ as a module over~$R$.  The next lemma
generalizes Section \ref{sub:braid} to show that, just as for the case
of the Burau cover~$\tinf{D}$, the first homology of~$\tinf{X}$ is a
\emph{free} module of finite rank over~$R$.
 
\begin{lemma}
  \label{dimlem}
  Suppose~$H_1(X)$ is a free Abelian group of finite rank~$r_1$, and
  let~$r_0$ be the number of connected components~$X_0 \subseteq X$
  such that~$\rho|_{H_1(X_0)}$ is not identically~$0$.
  Then~$H_1(\tinf{X})$ is a free~$R$-module of rank~$r = r_1 - r_0$.
\end{lemma}
\begin{proof}
  Let~$X' \subset X$ be a set consisting of one point in each
  connected component of~$X$, and let~$Z = (\tinf{p})^{-1}(X')$.
  Let~$\gamma_1, \dots, \gamma_{r_1}$ be a collection of loops in~$X$,
  each based at some point of~$X'$, whose homology classes
  generate~$H_1(X)$, and lift each~$\gamma_j$ to a
  path~$\tilde{\gamma}_j$ in~$\tinf{X}$.  Because~$\tinf{X}$ is a
  covering space with deck group~$\{T^i\}$, the homology classes of
  paths of the form~$T^i \tilde{\gamma}_j$ generate the relative
  homology group~$H_1(\tinf{X},Z)$ as a vector space.  In other words,
  the homology classes of the~$\tilde{\gamma}_j$'s
  generate~$H_1(\tinf{X},Z)$ as an~$R$-module.  Since~$H_1(X)$ is
  free, it follows that there can be no torsion in this~$R$-module,
  and therefore~$H_1(\tinf{X}, Z)$ is free and of rank~$r_1$.
  Moreover, since~$Z$ is discrete, the
  inclusion~$(\tinf{X},\emptyset)\subset(\tinf{X},Z)$ induces an
  injection~$H_1(\tinf{X}) \raw H_1(\tinf{X},Z)$ of~$R$-modules,
  so~$H_1(\tinf{X})$ is a free~$R$-module of rank at most~$r_1$.  The
  formula~$r = r_1 - r_0$ can now be proved by considering the long
  exact homology sequence of the pair~$(\tinf{X},Z)$. \QED\
\end{proof}

Let~$\zeta_1, \dots, \zeta_{r}$ be a basis for~$H_1(\tinf{X})$.  With
respect to this basis, the module isomorphism~$\tinf{h}_*$ is given by
a square matrix~$M = M(t) \in \GL(r,R)$ with entries in~$R$.  Note
that~$M$ depends on the morphism~$\rho$ used to define~$\tinf{X}$, on
the homotopy class of~$h$, and also on the choice of lift of~$h$
to~$\tinf{X}$.  If~$X=D$ is the~$n$-punctured disk,~$h$ is a
representative the braid~$\beta$ on~$n$ strings,~$\rho$ is the
homomorphism~$\tau$ of Section \ref{sub:braid} and~$\tinf{h}$ is the
preferred lift of~$h$ to~$\tinf{D}$, then we have~$M = B(\beta)$, the
reduced Burau matrix of~$\beta$.

For a complex number~$\eta \in \C$, denote by~$M(\eta)$ the complex
matrix obtained from~$M$ by substituting~$\eta$ in place of~$t$.
If~$\eta \neq 0$ then this matrix is invertible, just as~$M$ is, and
so it acts as a linear isomorphism of~$\C^{r}$ to itself.  We denote
by~$\RBSkC$ the complexification of $\RBSk$ and treat $\RBSkC$ as a
subspace of~$H_1(\tk{X},\C)$.

The next lemma connects the action on the $k$-fold cover to substitutions
of complex $k^{th}$ roots of unity into the matrix $M(t)$. It is based on
well-known, elementary facts. Depending on the chosen perspective,
it follows,  for example, from the splitting of a representation
of $Z_k$ into the sum of irreducibles, or from the invertibility of
the order-$k$ discrete Fourier transform. Rather than abstract the
necessary algebra and then apply it to the situation at hand, it 
is simpler to maintain an algebraic topology perspective and give
a direct proof using the invariance of an eigen-decomposition.

\begin{lemma}\label{eiglem}
  Let~$T$ be the generator of the deck group for the
  covering~$\tk{p}:\tk{X}\raw X$, and let~$\tk{h}$ and~$\tinf{h}$ be
  the lifts of~$h$ to~$\tk{X}$ and~$\tinf{X}$.  The eigenvalues
  of~$T_*$ restricted to ${\RBSkC}$ are~$1, \eta_k, \eta_k^2, \dots,
  \eta_k^{k-1}$ where $\eta_k = e^{2\pi i/k}$.  Denote by~$E_0, \dots,
  E_{k-1}$ the corresponding eigenspaces in~$\RBSkC$.  Then each
  subspace $E_{m}$ is $\tk{h}_*$-invariant, and the action of
  $\tk{h}_*$ on~$E_{m}$ is given by the matrix~$M(\eta_k^m)$, obtained
  by substituting~$\eta_k^m$ into the matrix~$M(t)$ of $\tinf{h}_*$.
\end{lemma}
\begin{proof}
  Let~$\Rk$ be the ring of all Laurent polynomials in a variable~$s$
  which satisfies~$s^k = 1$, so $\Rk$ is isomorphic to the group ring
  $\Z[\Z_k]$ .  As with ~$H_1(\tinf{X})$, we can treat~$H_1(\tk{X})$
  as a free module over~$\Rk$. Note that~$T_*$ and~$\tk{h}_*$ act by
  module isomorphisms, and $\RBSk$ is an $\Rk$-submodule that is
  invariant under both~$T_*$ and~$\tk{h}_*$.

  Since $T_*$ has order~$k$ its eigenvalues are as given. A general
  element of~$\RBSkC$ has the form~$\sum_{i \in \Z_k} s^i
  \left(\sum_{j=1}^{n-1} a_{i,j} \zeta_j^{(k)}\right)$, for complex
  numbers~$a_{i,j} \in \C$.  For $m = 0, \dots, k-1$, let~$E_m$ be the
  set of elements which when written in this form satisfy~$a_{i+1,j} =
  \eta_k^{-m} a_{i,j}$ for each~$i \in \Z_k$ and all $j$.  Since~$T_*$
  acts as an~$\RkC$-module homomorphism on $\RBSkC$ by multiplication
  by $s$, the set~$E_m$ consists of eigenvectors of~$T_*$ with
  eigenvalue~$\eta_k^m$.  Further, the dimension of~$E_m$ as a complex
  vector space is~$r$ and so $E_0\oplus\dots\oplus E_{k-1} = \RBSkC$.

  Since $\tk{h}_*$ commutes with~$T_*$, each~$E_m$ is
  $\tk{h}_*$-invariant.  By definition, the matrix~$M \in
  \GL(r, R)$ is the matrix of~$\tinf{h}_*: H_1(\tinf{X}) \raw
  H_1(\tinf{X})$ relative to the chosen basis~$\{\zeta_i\}$
  of~$H_1(\tinf{X})$.  We decompose this matrix as
  \begin{equation*}
    M = \sum_{i \in \Z} t^i M_i,
  \end{equation*}
  with each~$M_i\in\Mat(r, \Z)$.  Projecting this action to
  $\RBSkC$ we have that $\tk{h}_*$ acts as an $\RkC$-module
  homomorphism on $\RBSkC$ by the matrix
  \begin{equation}
    \label{eq:matrix of h_*}
    M(s) := \sum_{i \in \Z_k} s^i
    \left(\sum_{l \in \Z} M_{lk+i} \right) \in \GL(r,\Rk).
  \end{equation}
  On the other hand, since~$\eta_k^{k} = 1$, the matrix~$M(\eta_k^m)$
  is given by
  \begin{equation}\label{eqsub}
    M(\eta_k^m)
    = \sum_{i \in \Z} \eta_k^{mi} M_i \\
    = \sum_{i=0}^{k-1} \eta_k^{mi} \left( \sum_{l \in \Z} M_{lk+i} \right).
  \end{equation}
  If~$v \in E_m$ then by \eqref{eq:matrix of h_*} and \eqref{eqsub} we have
  \begin{align*}
    \tk{h}_* (v)
    & = \sum_{i \in \Z_k} s^i \left(\sum_{l \in \Z} M_{lk+i} \right) \cdot  v \\
    &= \sum_{i \in \Z_k} \left(\sum_{l \in \Z} M_{lk+i} \right) \cdot \eta_k^{mi} v \\
    &= M(\eta_k^m) \cdot v
  \end{align*}
  as claimed.
  \QED\
\end{proof}

Although~$H_1(\tk{X})$ is larger than~$\RBSk$, the next lemma indicates that all
of the \emph{growth} of~$\tk{h}_*$ on~$H_1(\tk{X})$ occurs in~$\RBSk$.
\begin{lemma}\label{growthlem}
  With notation as above, the eigenvalues of $\tk{h}_*$ acting on
  $H_1(\tk{X})$ are those of its restriction to $\RBSk$ together with
  some roots of unity.  In particular, the spectral radius
  of~$\tk{h}_*$ equals that of its restriction to $\RBSk$.
\end{lemma}
\begin{proof}
  We return to treating~$H_1(\tk{X})$ and~$\tk{S}$ as Abelian groups.
  Since~$\tk{S}$ is~$\tk{h}_*$-invariant, we may write~$\tk{h}_*$ as a
  matrix of the form
  \begin{equation*}
    \tk{h}_* =
    \left(
      \begin{matrix}
        A & B \\
        0 & C
      \end{matrix}
    \right)
  \end{equation*}
  in which~$A$ represents the action of~$\tk{h}_*$ on~$\tk{S}$,
  and~$C$ its action on~$\tk{h}_*$ on~$H_1(\tk{X}) / \tk{S}$.  We will
  construct a basis for~$H_1(\tk{X}) / \tk{S}$ whose elements are
  permuted by~$\tk{h}_*$.  In particular the eigenvalues of~$C$ are
  all roots of unity, so this will prove the lemma.

  To construct the desired basis, let~$Y$ be a connected component
  of~$\tk{X}$ and let~$X_0 = \tk{p}(Y)$.  The image of~$H_1(X_0)$
  under the morphism~$\rho:H_1(X)\raw\Z$ defining~$\tinf{X}$ is some
  subgroup~$a\Z$ of~$\Z$.  The chosen connected component~$Y$ then
  corresponds to a coset of~$a \Z_k$ in~$\Z_k$, and so the deck
  transformations which preserve~$Y$ are precisely those of the
  form~$T^i$ with~$i \in a \Z_k$.  

  Choose a loop~$\gamma$ in~$X_0$ such that~$\rho([\gamma]) = a$.
  This loop~$\gamma$ may lift in more than one way to a path in~$Y$.
  Let~$\tilde{\gamma}$ be one such lift.  We define an
  element~$\gamma(Y) \in H_1(\tk{X}) / \tk{S}$ by
  \begin{equation*}
    \gamma(Y) = \left[\sum_{i \in a \Z_k} T^j \tilde{\gamma} \right] + \tk{S}.
  \end{equation*}
  where, as usual, the square brackets denote taking the homology
  class.  Any other lift of~$\gamma$ to~$Y$ is of the form~$T^i
  \tilde{\gamma}$ for~$i \in a\Z_k$, so~$\gamma(Y)$ does not depend on
  the choice of lift of~$\gamma$.

  We claim that~$\gamma(Y)$ also does not depend on the chosen
  loop~$\gamma$.  For suppose that~$\gamma'$ is another loop in~$X_0$
  satisfying~$\rho([\gamma']) = a$.  Then~$\rho(\gamma - \gamma') = 0
  \in \Z$, so the closed~$1$-chain~$\gamma' - \gamma$ lifts to a
  closed~$1$-chain in~$\tinf{X}$, representing an element
  of~$H_1(\tinf{X})$.  If the lift~$\tilde{\gamma}'$ of~$\gamma'$ is
  chosen so that~$\tilde{\gamma}' - \tilde{\gamma}$ is also a
  closed~$1$-chain in~$\tk{X}$, then it follows that~$[\tilde{\gamma}'
  - \tilde{\gamma}] \in \tk{S}$.  Therefore also~$\sum_{i} T^i
  (\tilde{\gamma}' - \tilde{\gamma}) \in \tk{S}$, proving
  that~$\gamma'(Y) - \gamma(Y) = 0 \in H_1(\tk{X}) / \tk{S}$.

  Note that if~$a = 0$ then~$\gamma(Y) = 0 \in H_1(\tk{X}) / \tk{S}$,
  since in that case~$Y$ lifts to a component of~$\tinf{X}_0$
  homeomorphic to~$X_0$.  On the other hand if~$a \neq 0$ then the
  deck group of any connected component of~$\tinf{X}_0$ over~$X_0$ is
  infinite cyclic, so no lift of~$\sum_{i \in a\Z_k} T^i \tilde{\gamma}$ can
  represent a homology class in~$H_1(\tinf{X})$.  Hence~$\gamma(Y) =
  0$ iff~$a = 0$.  A dimension calculation similar to that in Lemma
  \ref{dimlem} now shows that the set~
  \begin{equation*}
    \{\gamma(Y)\st \text{$Y$ is a connected component of~$\tk{X}$}\}
  \end{equation*}
  generates~$H_1(\tk{X}) / \tk{S}$. To finish the proof, it is enough to note
  that, since~$\tk{h}$ permutes the connected components of~$\tk{X}$, and
  since~$\rho h_* = \rho$, the map~$\tk{h}_*$ must also permute the~$\gamma(Y)$'s.
  \QED
\end{proof}

Lemma~\ref{growthlem} together with Lemma~\ref{eiglem} yields
\begin{theorem}\label{Main Algebra Theorem}
  Let~$h:X\raw X$ be a homeomorphism of the locally path-connected,
  semi-locally simply connected topological space~$X$, whose first
  homology group we assume to be free and of finite rank.
  Suppose~$\rho:H_1(X) \raw \Z$ is a homomorphism which
  satisfies~$\rho h_* = \rho$, and let~$\tinf{X}$ and~$\tk{X} =
  \tinf{X} / T^k$ denote the covering spaces over~$X$ corresponding
  to~$\rho$ and~$\xi_k\circ \rho$, with covering
  projection~$\tk{q}:\tinf{X}\raw\tk{X}$.  Let~$\tinf{h}$ and~$\tk{h}$
  denote lifts of~$h$ to these covering spaces.  If~$M = M(t) \in
  GL(r,R)$ denotes the matrix of~$\tinf{h}_*:H_1(\tinf{X})\raw
  H_1(\tinf{X})$ as an~$R$-module isomorphism, then the action
  of~$\tk{h}_*$ on the invariant subspace~$\tk{S}_{\C} =
  \tk{q}_*(H_1(\tinf{X},\C))$ is given by the direct sum
  \begin{equation*}
    \tk{h}_* = M(1) \oplus M(\eta_k) \oplus \cdots \oplus M(\eta_k^{k-1})
  \end{equation*}
  where~$M(\eta_k^j)$ denotes the complex matrix obtained by
  substituting~$\eta_k^j = e^{2 \pi i j/k}$ into~$M$.  Furthermore,
  any eigenvector of~$\tk{h}_*$ not lying in~$\tk{S}$ has eigenvalue
  which is a root of unity.
\end{theorem}

\section{Entropy and  first homology}\label{sect:ent}
The topological entropy, $\htop(f)$, is a well-known measure of the
complexity of the dynamics of a self-map $f$ of a compact metric
space.  See \cite{AdlerKonheimMcAndrew}, \cite{Denker1976} or \cite{Katok1995} for more
information.  The next well-known lemma contains part of the main idea
in Fried's paper \cite{Fried1986}: passing to a finite cover often
allows one to detect more growth on homology.

\begin{lemma} \label{toplem} If $f:Y\raw Y$ is a continuous map of the
  compact manifold~$Y$, then~$$\htop(f) \geq \sup \left\{
    \log\left(\spec \left(\tf_*\right) \right) \right\}$$
  where the
  supremum is over all lifts, $\tf$, of $f$ to a finite
  cover~$\tilde{Y}$ and~$\tf_*$ is the action of $\tilde{f}$ on first
  homology of the cover~$H_1(\tilde{Y};\R)$.
\end{lemma}

\begin{proof}
  The lemma follows directly from two classic results.  Manning proved
  in \cite{Manning1974} that $h_{top}(f) \geq \log(\spec(f_*))$ for
  any continuous $f:Y\raw Y$ and Bowen proved in \cite{Bowen1971} that
  entropy is preserved under finite to one factors, in particular,
  $\htop(f) = \htop(\tf)$ for any lift $\tf$ of $f$ to a finite
  cover. \QED\
\end{proof}

For a braid $\beta$ we define $\htop(\beta) = \inf\{\htop(h): h\in
\beta\}$.  Now Theorem \ref{Main Algebra Theorem} says that the
maximal spectral radius of $ B(e^{2 \pi i j/k})$ for $0\leq j < k$
gives the spectral radius on the first homology of the lift to the
$k$-fold Burau cover $\tk{D}$ of an $h\in \beta$.  Thus by
Lemma~\ref{toplem} we have what was referred to in  the Introduction
as the \de{Burau estimate}.

\begin{equation*}
  \htop(\beta) \geq
  \sup \left\{\log(\spec(B(e^{2 \pi i j/k}))) : j/k \in\rationals \right\}.
\end{equation*}
Since the entries in $B(t)$ are polynomials, $\spec(B(t))$ is
continuous in $t$. Thus we have the result of Fried
\cite{Fried1986} and Kolev \cite{Kolev1989}:
\begin{lemma}
  If $\beta\in B_n$ with Burau representation
  $B(t)$, then 
  \begin{equation}\label{burest}
    \htop(\beta) \geq \sup\{\log(\spec(B(\eta))) : \eta\in\circle\}.
  \end{equation}
\end{lemma}

A fundamental result in
Nielsen-Thurston theory says that $\htop(\beta) = \htop(\Phi)$ where
$\Phi$ as given in Theorem~\ref{thurstons theorem} is the Thurston
representative in the isotopy class $\beta$.  A natural question is
whether this value is detected by the action on homology in some
finite cover, or in the current context, whether the estimate in
\eqref{burest} is ever sharp.

To investigate this question we must first understand a simpler
question, namely, for a \pA\ map $\phi$ of a surface $X$, when is
$\htop(\phi) = \spec(\phi_*)$?  In other words, when is Manning's
estimate in \cite{Manning1974} sharp for a \pA\ map? The answer turns
out to depend exactly on the orientability of the $\phi$-invariant
foliations.  It is well-known that the oriented, measured foliation
$\cF^u$ gives rise to a homology class $v\in H_1(X;\R)$, for example,
 as an asymptotic direction as in Schwartzman (\cite{Schwartzman1957}) or a
geometric current as in Ruelle and Sullivan(\cite{RuelleSullivan}).
Briefly, $v$ is the average direction in homology obtained from
flowing along the one-dimensional leaves of $\cF^u$.  This average
exists and is unique because the unstable foliation of a \pA\ map is
uniquely ergodic (\cite{Orsay1979}).

Now if the unstable foliation $\cF^u$ gives rise to $v^u\in H_1(X;\R)$, then
since $\phi_* \cF^u = \lambda \cF^u$, on first homology we have $\phi_* v^u =
\lambda v^u$.  Thus $v^u$ is an eigenvector of $\phi_*$ with eigenvalue
$\lambda$, and so $\spec(\phi_*) \geq \lambda$. On the other hand, since \pA\
maps have entropy equal to the logarithm of their expansion constant,
Lemma~\ref{toplem} yields that $\lambda \geq \spec(\phi_*)$.  Thus when the
invariant foliations are orientable we see that the spectral radius on first
homology gives the entropy of a \pA\ map.  The converse of this fact doesn't
seem as well-known so we include a proof below for completeness.

\begin{lemma}
  \label{lemma: orientability vs eigenvalues}
  Suppose~$\phi:X\raw X$ is a \pA\ \homeo\ of the orientable
  surface~$X$ having~$\ell$ connected components (which, according to
  the definition of pseudo-Anosov map, must be permuted cyclically
  by~$\phi$), and let~$\lambda$ be the expansion constant of~$\phi$.
  \begin{tightlist}
  \item{(a)} The pA map~$\phi$ has orientable invariant foliations in some
    (hence each) connected component if and only if $\spec(\phi_*) = \lambda$.
  \item{(b)} Suppose~$\phi$ has oriented invariant foliations, and let~$\epsilon
    = \pm 1$ according to whether~$\phi$ preserves or reverses the orientation
    of the unstable foliation.  Then each complex number of the form~$\epsilon
    e^{2 \pi i j/\ell} \lambda$ is a simple eigenvalue of~$\phi_*$, and every
    other eigenvalue~$\mu$ satisfies~$|\mu| < \lambda$.
  \end{tightlist}
\end{lemma}

\begin{proof}
  We first suppose that~$X$ is connected, i.e.\ that~$\ell=1$, and prove
(b) and then (a).  The proof when~$\ell > 1$ is an easy modification.

 When $\ell = 1$, $(b)$
  states that if the invariant foliations of~$\phi$ are oriented, then~$\epsilon
  \lambda$ is a simple eigenvalue of~$\phi_*$, and every other eigenvalue is
  smaller in modulus.  Because~$X$ is orientable, orientability of the unstable
  foliation of~$\phi$ is equivalent to orientability of the stable foliation
  of~$\phi$ (\cite{CamachoAndNeto1985}), so we need only talk about
  orientability of the former.

  In \cite{Rykken1999} Rykken shows that if the \pA\ map~$\phi$ has
  an orientable unstable foliation, then except for zeros and roots of
  unity, the eigenvalues of its action on first homology are the same
  as those of its Markov transition matrix, $A_{\mathcal{P}}$, with
  respect to any Markov partition~$\mathcal{P}$. Using this result (b)
  then follows from the Perron-Frobenius Theorem applied to
  $A_{\mathcal{P}}$ (see, for example, \cite{Kitchens1998}).

  We noted above the theorem how one implication in (a) follows from
  treating the foliation as an asymptotic cycle. We prove the
  contrapositive of the converse and so assume that the unstable
  foliation~$\cF^u$ of $\phi$ is not orientable.
  Let~$\overrightarrow{X} \stackrel{\op}{\raw} X$ be the orientation
  double-cover of~$\cF^u$ constructed in Section~\ref{sub:branch};
  $\phi$ lifts to a pseudo-Anosov $\orient{\phi}$
  of~$\overrightarrow{X}$ with unstable foliation
  denoted~$\overrightarrow{\cF}^u$.  The leaves of
  $\overrightarrow{\cF}^u$ of necessity project under~$\op$ to those
  of~$\cF^u$.  Since $\orient{\phi}$ also has expansion
  constant~$\lambda$, by part (b) $\ophistar:H_1(\overrightarrow{X})
  \raw H_1(\overrightarrow{X})$ has a simple eigenvalue equal
  to~$\epsilon \lambda$ and all its other eigenvalues have smaller
  modulus.

  If we let $\theta$ denote the generator of the deck group of $\orient{X}$,
  then $\theta$ will take leaves of $\orient{\cF^u}$ to leaves reversing the
  orientation while preserving the transverse measure, and so $\theta_* v^u =
  -v^u$.  Now since~$\theta^2 = \id$, $H_1(\overrightarrow{X}, \R) = V_+ \oplus
  V_-$ with $V_\pm$ the eigenspace of $\theta_*$ with eigenvalue $\pm 1$.
  Moreover, we now show that $V_-$ is precisely the kernel of the map induced on
  homology by the covering projection $\opstar:H_1(\overrightarrow{X},\R)\raw
  H_1(X,\R)$.  Suppose~$\tau$ is a~$1$-cycle in~$\overrightarrow{X}$ avoiding
  the singular points of~$\overrightarrow{\cF}$ such that~$\opstar [\tau] = 0
  \in H_1(X)$.  Then we may write~$\opstar \tau = \partial \Delta$ where~$\Delta
  = \sum_i a_i \Delta_i$ is some $2$-chain.  Lift each~$\Delta_i$ to a
  $2$-simplex~$\Delta'_i$ in~$\overrightarrow{X}$ and
  let~$\overrightarrow{\Delta} = \sum_i a_i \left(\Delta'_i + \theta
    \Delta'_i\right)$.  Clearly,~$\partial \overrightarrow{\Delta} = \tau +
  \theta \tau$ proving that~$[\tau] + [\theta \tau] = 0 \in
  H_1(\overrightarrow{X})$.  In other words~$[\tau] \in V_-$.  So~ $\ker \opstar
  \subset V_-$.  Conversely if~$v \in V_-$, then~$2 \opstar(v) = \opstar(v) +
  \opstar(\theta(v)) = \opstar(v) - \opstar(v) = 0$, and hence~$V_- = \ker
  \opstar$ as claimed.

  Thus~$\phi_*$ acting on~$H_1(X,\R)$ is conjugate to~$\ophistar$
  acting on~$V_+$.  Since~$v^u \in V_-$ and the eigenvalue~$\epsilon
  \lambda$ is geometrically simple, it follows that every eigenvalue
  of~$\ophistar$ having an eigenvector lying in~$V_+$
  has modulus strictly smaller than~$\lambda$, so that~$\spec(\phi_*)
  < \lambda$ as required.  \QED
\end{proof}

\section{The Burau estimate for pseudo-Anosov braids}\label{sect:burest}
A braid~$\beta \in B_n$ is said to be \emph{pseudo-Anosov} if in the
Thurston normal form $\Phi\in\beta$ of Theorem \ref{thurstons theorem}
the set~$A$ consists of just one annulus~$a$ which is the collar
of~$\bd D$, and the restriction of~$\Phi$ to the complement of~$a$ is
a pseudo-Anosov map $\phi$.  In this case we will call~$\Phi$ a
\emph{collared pseudo-Anosov} map and consider its invariant
foliations to be those of $\phi$.

For simplicity of notation we let $r(\theta) = \spec(B(\exp{2 \pi i
  \theta}))$ and just consider $\theta\in[0, 2 \pi)$.  If $\beta$ is
\pA\ with expansion constant $\lambda = e^{\htop(\Phi)}$,
Lemma~\ref{toplem} shows that $r(\theta) \leq \lambda$ for all
$\theta\in[0, 2 \pi)$.  In this section we investigate when the Burau
estimate \eqref{burest} is sharp for a \pA\ braid and the substitution
of a root of unity, that is, when do we have $r(j/k) = \lambda$ for
some $j/k\in\rationals$. The first observation is that according to
Theorem~\ref{Main Algebra Theorem} when $r(j/k) = \lambda$ occurs,
then the action of ${h}^{(k)}_*$ on the first homology of the $k$-fold
Burau cover $H_1(\tk{D})$ has some eigenvalue of modulus~$\lambda$. By
Lemma~\ref{lemma: orientability vs eigenvalues}, this eigenvalue must
be either~$\lambda$ or~$-\lambda$, it must be simple, and all other
eigenvalues of ${h}^{(k)}_*$ must be of smaller modulus.  So, in
particular, $r(\ell/k)< \lambda$ for all $\ell \not = j$.

The second observation is that the function $r(\theta)$ is periodic
with period one and in addition, since the coefficients of the 
Laurent polynomial entries of $B(t)$ are real, $\spec(B(\eta)) =
\spec(B(\overline{\eta}))$, where~$\overline{\eta}$ denotes the complex
conjugate of~$\eta$. Thus $r$ is an even, one-periodic function and so
it is even about $1/2$.

Putting the two observations together we have that that for a given
$k$, $r(j/k) = \lambda$ for at most one $j$ and since $r(1/2 +x) =
r(1/2 -x)$, that can only happen if $j/k = 1/2$. Thus if $k$ is even,
the only possibility for the Burau estimate to be sharp is that $r(
(k/2) / k) = \lambda$ and so $r( j / k) < \lambda$ for other $j$. On
the other hand, if $k$ odd, $r( j / k) < \lambda$ for all $j$, so the
Burau estimate is never sharp for substitutions with $k$ odd.

In the next Proposition we connect these observations with
the structure of the invariant foliations from 
Section~\ref{sub:branch} and 
Lemma~\ref{lemma: orientability vs eigenvalues}.

\begin{theorem} \label{propn: pA burau estimate} 
  Suppose~$\beta \in B_n$ is a braid represented by the
  collared pseudo-Anosov map $\Phi:D \raw D$ having
  expansion factor~$\lambda > 1$.  The following are
  equivalent:
  \begin{tightlist}
  \item{(a)} $\spec\left(B(e^{2\pi i j/k})\right) = \lambda$
    for some~$k > 0$ and some~$0 \leq j < k$;
  \item{(b)} $\spec\left(B(-1)\right) = \lambda$ 
    and~$-1$ is the only root of unity for which this occurs.
  \item{(c)} The invariant foliations~$\cF^u$ and~$\cF^s$ have an
    odd-order singularity at each puncture of~$D$, and all singularities
    of~$\cF^u$ and~$\cF^s$ in the interior of~$D$ have even order.
  \item{(d)} $D^{(2)}$ is the orientation double-cover
    of~$\cF^u$ and~$\cF^s$ (after removing the collaring).
  \end{tightlist}
\end{theorem}

\begin{proof}
  The observations above the Theorem
  prove the equivalence of (a) and (b).

  Let $\tPhit$ be the preferred lift of~$\Phi$ to~$D^{(2)}$ and
  $\tPhit_*$ its action on $H_1(D^{(2)})$. By Theorem~\ref{Main Algebra
    Theorem}, the eigenvalues of~$\tPhit_*$ are those of~$B(1) \oplus
  B(-1)$, together with some roots of unity.  Since~$B(1)$ is a
  permutation,~$\lambda$ (or~$-\lambda$) is an eigenvalue of $\tPhit_*$
  if and only if it is an eigenvalue of~$B(-1)$. In addition, by
  Lemma~\ref{lemma: orientability vs eigenvalues}, since $\tPhit$ is
  collared \pA, one of~$\lambda$ or~$-\lambda$ is an eigenvalue of
  $\tPhit_*$ if and only the invariant foliations of~$\tPhit$ are
  oriented. Thus (b) holds exactly when the lift of $\cF^u$ to $D^{(2)}$
  is oriented, and so in particular (d) implies (b).

  As remarked above Lemma~\ref{orientlemma}, if $\gamma$ is a small
  loop surrounding a singularity $P$, then recalling that
  $\orient{\rho}$ denotes the epimorphism defining the orientation
  double cover, we have $\orient{\rho}([\gamma]) = 1$ iff $P$ has odd
  order.  On the other hand, the epimorphism $\tau_2:H_1(D)\raw \Z_2$
  which yields the two-fold Burau cover $D^{(2)}$ is defined by
  $\tau_2([\gamma']) = 1$ for $\gamma'$ a small loop around a
  puncture. This shows the equivalence of (c) and (d).

  Now we will show (b) implies (c) by proving its contrapositive.  If
  $\cF^u$ has an interior singularity of odd order, by
  Lemma~\ref{orientlemma} the lifted foliations to $D^{(2)}$ are not
  oriented. Now note that by Euler-Poincar\'{e}-Hopf formula
  \eqref{EPHFormula}, $\cF^u$ always has a one-pronged singularity at
  some puncture~$x_i$ of $D$.  If $\cF^u$ has an even order
  singularity at some other puncture~$x_j$, we consider the homotopy
  class~$\alpha = \gamma_j\I \gamma_i$. We then have
  $\orient{\rho}(\alpha) = 0 + 1 = 1$ and $\tau_2 (\alpha) = -1 + 1 =
  0$. Thus there can be no homomorphism $\delta$ with
  $\delta\circ\rho_2 = \orient{\rho}$, and so again by
  Lemma~\ref{orientlemma}, the lifted foliations to $D^{(2)}$ are not
  oriented.  Since we have just seen that (b) holds exactly when the
  lift of $\cF^u$ to $D^{(2)}$ is oriented, we have that (b) implies
  (c).  \QED\
\end{proof}

\begin{remark}\label{rk:song}
Lemma 5 and Theorem 7 in \cite{Song2002} contain a portion of
the above Theorem in slightly different language, namely, the implications
$(c) \implies (d)$ and $(d) \implies (b)$.
\end{remark}
 
\section{The Burau estimate for reducible braids} \label{sect:reducible}

\subsection{Pseudo-Anosov components and Burau orientability}

The Thurston normal form~$\red = \red_1 \cup \dots \cup \red_m$ of a
general braid~$\beta \in B_n$ can be quite complicated.  Using the
notation of Theorem \ref{thurstons theorem}, there may be several
reducing annuli forming the set~$A$, each
component~$\red_i=\red|_{X_i}$ of~$\red$ may be either periodic or
pseudo-Anosov, and~$X_i$ itself may have one or more connected
components.  In this section we will consider this more complicated
decomposition and obtain general results on the sharpness of the
entropy bound provided by the Burau representation.

The simplest case is when each~$\red_i$ is periodic, and
so~$\htop(\red) = 0$. Since all the eigenvalues of~$\red_*$ are roots
of unity, we have~$\spec(B(1)) = 1$, giving a sharp bound on the
entropy.  Thus from now on we suppose that some component of~$\red$ is
pseudo-Anosov, and we let~$\lambda = e^{\htop(\red)}$.  The main
theorem of this section gives necessary and sufficient conditions for
a root of unity~$\omega$ to make the Burau estimate~\eqref{burest}
sharp, i.e.\ to satisfy~$\spec(B(\omega)) = \lambda$.  

To state these conditions, we first define for each component~$\red_i$
a positive integer~$a_i$ expressing the manner in which the
surface~$X_i$ on which~$\red_i$ is supported surrounds the punctures
of~$D$. Begin by choosing a component~$\red_i$ of~$\red$.  The
supporting surface~$X_i$ of~$\red_i$ may have several connected
components; by definition, these are permuted cyclically by~$\red_i$.
Let~$X_{i0}$ be one of these components.  Thus~$X_{i0}$ is a finitely
punctured subdisk of~$D$ from which a finite collection of open
punctured subdisks have been deleted (see Figure
\ref{fig:areduciblebraid}(a)).  Let~$x_1, \dots, x_{r'}$ denote the
punctures in~$X_{i0}$, where~$0 \leq r' \leq n$, and let~$O_{r'+1},
\dots, O_r$ be the deleted subdisks.  For~$r' < j \leq r$ we
write~$m_j$ for the number of punctures of~$D$ which lie in~$O_j$; and
for~$1 \leq j \leq r'$ we define~$m_j = 1$.  Finally, we define~$a_i$
by
\begin{equation}
  \label{eq:a_i}
  a_i = \gcd(m_1, m_2, \dots, m_{r}).
\end{equation}
Because~$\red$ permutes the connected components of~$X_i$ cyclically,
and sends punctures to punctures, this definition is independent of
the choice of~$X_{i0}$.

The next lemma gives another interpretation of the number~$a_i$ which
follows from the definition of the morphism~$\tau:H_1(D)\raw
\Z$ associated to the covering space~$\tinf{D}$.
\begin{lemma}
  \label{lem:a_i}
  Let~$\ell_i$ be the number of connected components of the supporting
  surface~$X_i$ of~$\red_i$, and let~$X_{i0}$ be the chosen connected component
  of~$X_i$ as above.  Write~$\tinf{X}_i = (\tinf{p})^{-1}(X_i)$
  and~$\tinf{X}_{{i0}} = (\tinf{p})^{-1}(X_{i0})$.  The number~$a_i$ just defined
  is equal to the number of connected components of~$\tinf{X}_{{i0}}$.  Hence the
  number of connected components of~$\tinf{X}_i$ is~$\ell_i a_i$.
\end{lemma}
\begin{proof}
  According to Thurston's Theorem \ref{thurstons theorem}, the waist
  curve of any reducing annulus~$a \in A$ is neither null-homotopic
  nor homotopic to a puncture.  It follows that the inclusion~$X_i
  \subset D$ induces an injection~$H_1(X_i) \raw H_1(D)$, and so we can
  consider~$H_1(X_i)$ as a subgroup of~$H_1(D)$.  Then~$\tinf{X}_i
  \raw X_i$ is (isomorphic to) the covering space defined by the
  homomorphism~$\tau|_{H_1(X_i)}$.  By definition~$\tau$ sends a small
  clockwise loop around the puncture~$x_j$ (respectively, the hole~$O_j$)
  of~$X_{i0}$ to~$m_j \in \Z$, and therefore~$a_i$ is the generator
  of~$\im \tau|_{H_1(X_{i0})}$.  This proves the first statement of
  the lemma.  Since~$h$ permutes the connected components of~$X_i$
  cyclically and since~$\tau \red_* = \tau$, we have that $\tinf{X}_i$
  has the same number of connected components above each component
  of~$X_i$, proving the second statement. \QED
\end{proof}

Let us say that a pseudo-Anosov component~$\red_i$ of~$\red$ is
\emph{Burau orientable} provided the lifts of its invariant foliations
to some~$\tk{D}$ are orientable.  The main theorem of this section is:

\begin{theorem}
  \label{theorem:spectral radius for reducible braids}
  Let~$\red = \red_1 \cup \dots \cup \red_m$ be the Thurston normal form of the
  braid~$\beta$, and for each~$i$ such that~$\red_i$ is pseudo-Anosov, let~$a_i$
  be as in \eqref{eq:a_i}.  
  \begin{enumerate}
  \item The pseudo-Anosov component~$\red_i$ is Burau orientable if
    and only if (with notation as in the definition of~$a_i$) the
    invariant foliations of~$\red_i$ have a singularity of odd order
    at each puncture~$x_j$ of~$X_{i0}$, and a singularity of odd
    (even) order on the boundary of each deleted disk~$O_j$ such
    that~$m_j / a_i$ is odd (even).
  \item Let~$I$ be the set of~$1 \leq i \leq m$ such that~$\red_i$ is
    pseudo-Anosov, Burau orientable, and satisfies~$\htop(\red_i) =
    \htop(\red)$.  Then the set
      of roots of unity~$\omega$ for which~$\spec(B(\omega)) = \lambda$, is equal to the union
      over~$i \in I$ of the set of all~$a_i$-th roots of~$-1$.

      In particular, if~$I$ is empty, then~$\spec(B(\omega)) <
      \lambda$ for every root of unity~$\omega$.
    \end{enumerate}
\end{theorem}

Before embarking on the proof of Theorem~\ref{theorem:spectral radius
  for reducible braids}, we will illustrate the theorem with some
examples.

\subsection{Examples}

It is convenient to use the following notation.  Let~$n > 0$.
If~$i,n_1,n_2$ are positive integers with~$i+n_1+n_2-1 \leq n$, we
let~$\sigma_{i,n_1,n_2} \in B_n$ denote the braid which moves the
group of of~$n_1$ consecutive strings starting at string~$i$ over the
group of~$n_2$ consecutive strings starting at string~$i+n_1$:
\begin{equation*}
  \sigma_{i,n_1,n_2} =
  (\sigma_{i+n_1-1} \cdots \sigma_{i+n_1+n_2-2})
  (\sigma_{i+n_1-2} \cdots \sigma_{i+n_1+n_2-3})
  \dots
  (\sigma_{i} \cdots \sigma_{i+n_2-1}).
\end{equation*}
In particular~$\sigma_{i,1,1} = \sigma_i$ for all~$1 \leq i \leq n-1$.
\begin{example}
  \label{ex:1}
  Let~$n' \geq 1$, and define a braid~$\beta_{n'}$ on~$3n'$ strings by
  setting~
  \begin{equation*}
    \beta_{n'} = \sigma_{1,n',n'} \sigma_{n'+1,n',n'}^{-1}.
  \end{equation*}
  Figure \ref{fig:areduciblebraid}(b) shows~$\beta_{n'}$ when~$n'=3$.
\begin{figure}
  \begin{equation*}
    \psfrag{a}{(a)}
    \psfrag{b}{(b)}
    \psfrag{A}{$\sigma_{1,3,3}$}
    \psfrag{B}{$\sigma^{-1}_{4,3,3}$}
    \psfrag{h}{$\red_i$}
    \psfrag{x}{$x_1$}
    \psfrag{y}{$x_2$}
    \psfrag{X}{$X_{i0}$}
    \psfrag{Y}{$X_i$}
    \psfrag{Z}{$\red_i(X_{i0})$}
    \psfrag{O}{$O_{3}$}
    \scalebox{0.7}{
      \includegraphics{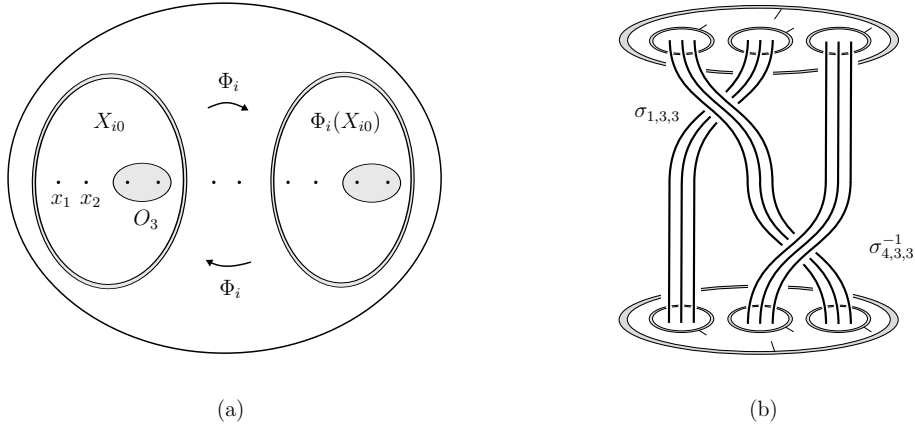}
    }
  \end{equation*}
  \caption{(a) A typical supporting surface~$X_i$.
    Here~$r'=2$,~$r=3$,~$m_3=2$ and~$a_i=1$.  (b) the
    braid~$\beta_{3}$ and its Thurston decomposition.}
  \label{fig:areduciblebraid}
\end{figure}
The Thurston normal form~$\red$ of~$\beta_{n'}$ reduces along~$4$
annuli~$A = A_0 \cup A_1 \cup A_2 \cup A_3$: there is one such annulus
collaring the boundary of~$D$, and one surrounding each of the three
groups of~$n'$ punctures.  The component of~$\red$ in the outer
connected component of~$D \setminus A$ is pseudo-Anosov.  Call this
component~$\red_1$ and its supporting surface~$X_1$.  The other
component of~$\red$ is periodic, and it cyclically permutes the three
inner connected components of~$D \setminus A$.  The
Euler-Poincar\'{e}-Hopf formula \eqref{EPHFormula} shows that the
invariant foliations of~$\red_1$ have four~$1$-pronged singularities,
one at each of the four boundary components of~$X_1$, and no other
singularities.  According to \eqref{eq:a_i} we
have~$a_1=m_1=m_2=m_3=n'$, and hence~$\red_1$ is Burau orientable by
the first part of Theorem \ref{theorem:spectral radius for reducible
  braids}.  

The growth rate of~$\beta_{n'}$ is~$\lambda = e^{\htop(\red_1)} \sim
2.618$.  Figure \ref{fig:specrad} shows the graph of the map
sending~$\theta \in [0,1]$ to the spectral radius of the substituted
Burau matrix~$B(e^{2 \pi i \theta})$ of~$\beta_{n'}$, in the two
cases~$n' = 8 = 2^3$ and~$n'=5$.

Here is one part of the justification of Theorem \ref{theorem:spectral
  radius for reducible braids} for this example.  According to Lemma
\ref{lemma: reducible vs orientability} below, if~$k$ is a multiple
of~$2n'$, the invariant foliations of~$\red_1$ lift to orientable
foliations in the~$k$-fold covering space~$\tk{X}_1$ of~$X_1$.
Moreover this covering space has exactly~$n'$ connected components,
each fixed by the lift of~$\red$.  Therefore by Theorem \ref{Main
  Algebra Theorem} and Lemma \ref{lemma: orientability vs eigenvalues}
we expect the matrices~$B(\omega)$, for~$\omega$ a~$k$-th root of
unity, to contribute exactly~$n'$ eigenvalues equal to~$\lambda$
or~$-\lambda$ (counted with multiplicity).  Theorem
\ref{theorem:spectral radius for reducible braids} states, in
addition, that it is precisely the~$n'$-th roots of~$-1$ which
contribute these eigenvalues, a fact clearly reflected in Figure
\ref{fig:specrad}.

\begin{figure}
  \begin{align*}
    \begin{array}{cc}
      \scalebox{0.45}{
        \includegraphics[bb=90 250 520 578,clip]{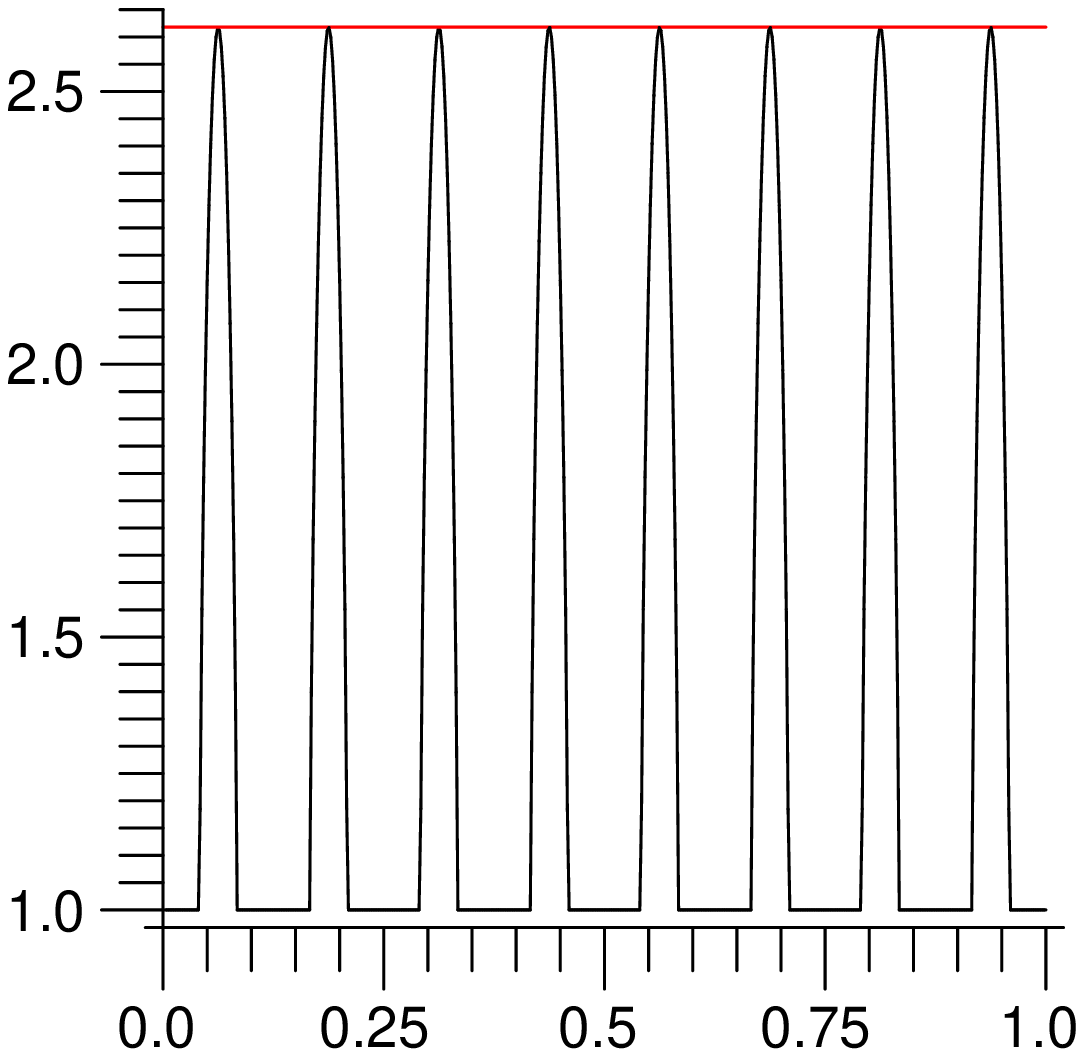}
      }
      &
      \scalebox{0.45}{
        \includegraphics[bb=90 252 520 578,clip]{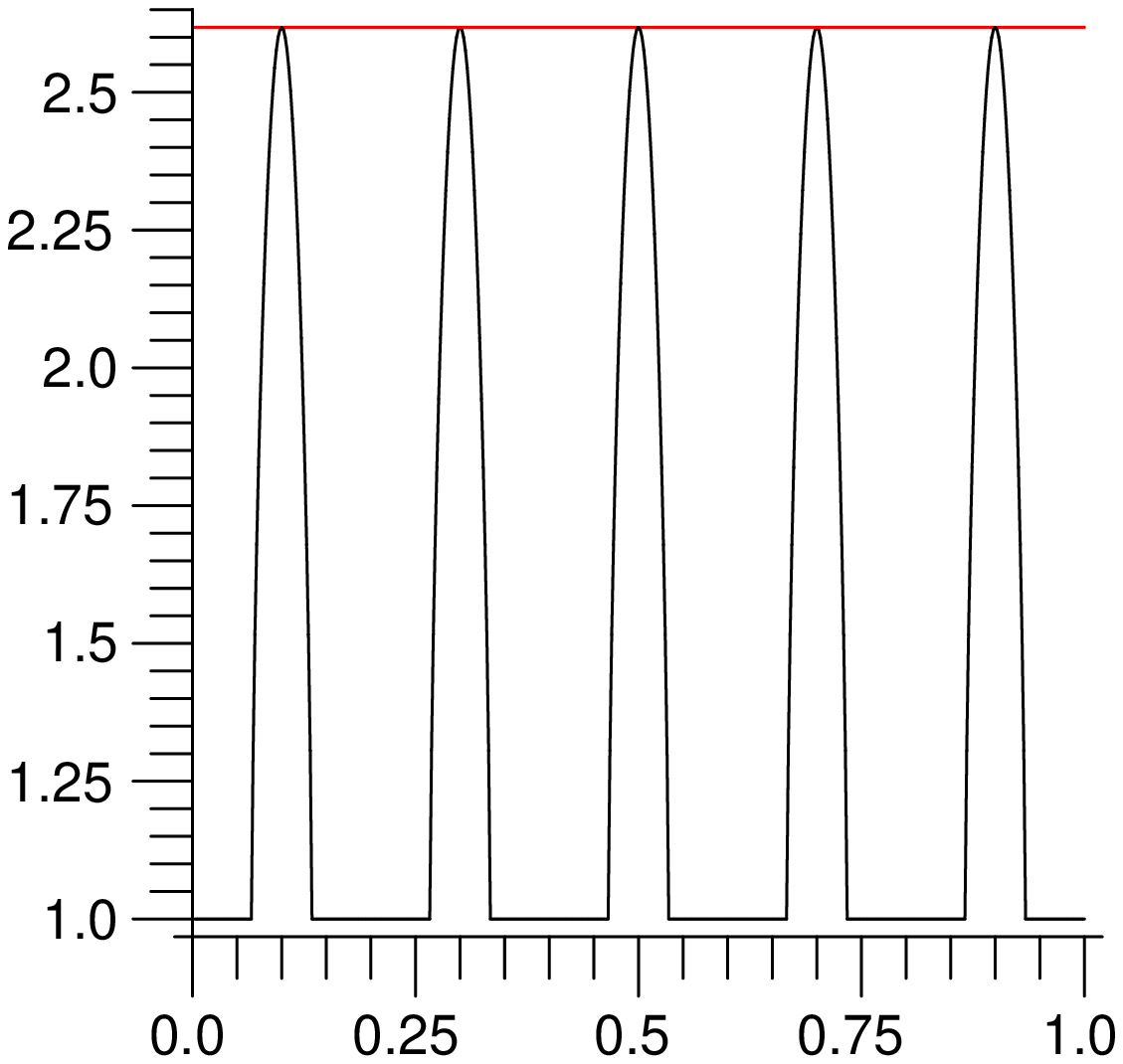}
      } \\
      (a) & (b)
    \end{array}
    \end{align*}
    \caption{The Burau estimate (a) for $\beta_8$ and (b)
      for~$\beta_5$.  The horizontal lines represent the growth
      rate~$\lambda \sim 2.618$.}
  \label{fig:specrad}
\end{figure}
\end{example}

\begin{example}
  Let~$\beta'_1 \in B_9$ and~$\beta'_2 \in B_8$ be the braids
  \begin{align*}
    \beta'_1 &= \sigma_{1,3,3} \cdot \sigma_{4,3,3}^2 \cdot \sigma_{1,3,3}^3 \\
    \intertext{and} 
    \beta'_2 &= \sigma_{1,3,3} \cdot \sigma_{4,3,2} \cdot \sigma_{4,2,3} \cdot \sigma_{1,3,3}^3.
 \end{align*}
 Note that~$\beta'_1$ respects the grouping of the punctures into
 consecutive groups of three, while~$\beta'_2$ respects the grouping
 of the punctures into two consecutive groups of three and one of two.
 As for the previous example, the Thurston normal form~$\red$
 of~$\beta'_i$ has one reducing annulus around each of the three
 groups of punctures.  Each such annulus encloses a periodic (in fact
 fixed) component of~$\red$.  If~$X_1$ denotes the complement of the
 reducing annuli, then~$\red_1 = \red|_{X_1}$ is pseudo-Anosov with
 growth rate~$\lambda \sim 5.828$.  Again the invariant foliations
 of~$\red_1$ have a one-pronged singularity on each boundary component
 of~$X_1$.

 Now for~$\beta'_1$ we have~$a_1=m_1=m_2=m_3=3$, and so Theorem
 \ref{theorem:spectral radius for reducible braids} shows
 that~$\red_1$ is Burau orientable.  Therefore the Burau estimate is
 sharp at all of the cubic roots of~$-1$ (see Figure
 \ref{fig:specrad2}(a)). For~$\beta'_2$, however, we have~$m_1=m_2=3$
 and~$m_3=2$, and so~$a_1 = \gcd(2,3) = 1$.  Thus the conditions of
 the first part of Theorem \ref{theorem:spectral radius for reducible
   braids} fail to hold and~$\red_1$ is not Burau orientable, so the
 Burau estimate is never sharp (see Figure \ref{fig:specrad2}(b)).
\begin{figure}
  \begin{equation*}
    \begin{array}{cc}
      \scalebox{0.45}{
        \includegraphics[bb=80 250 520 575,clip]{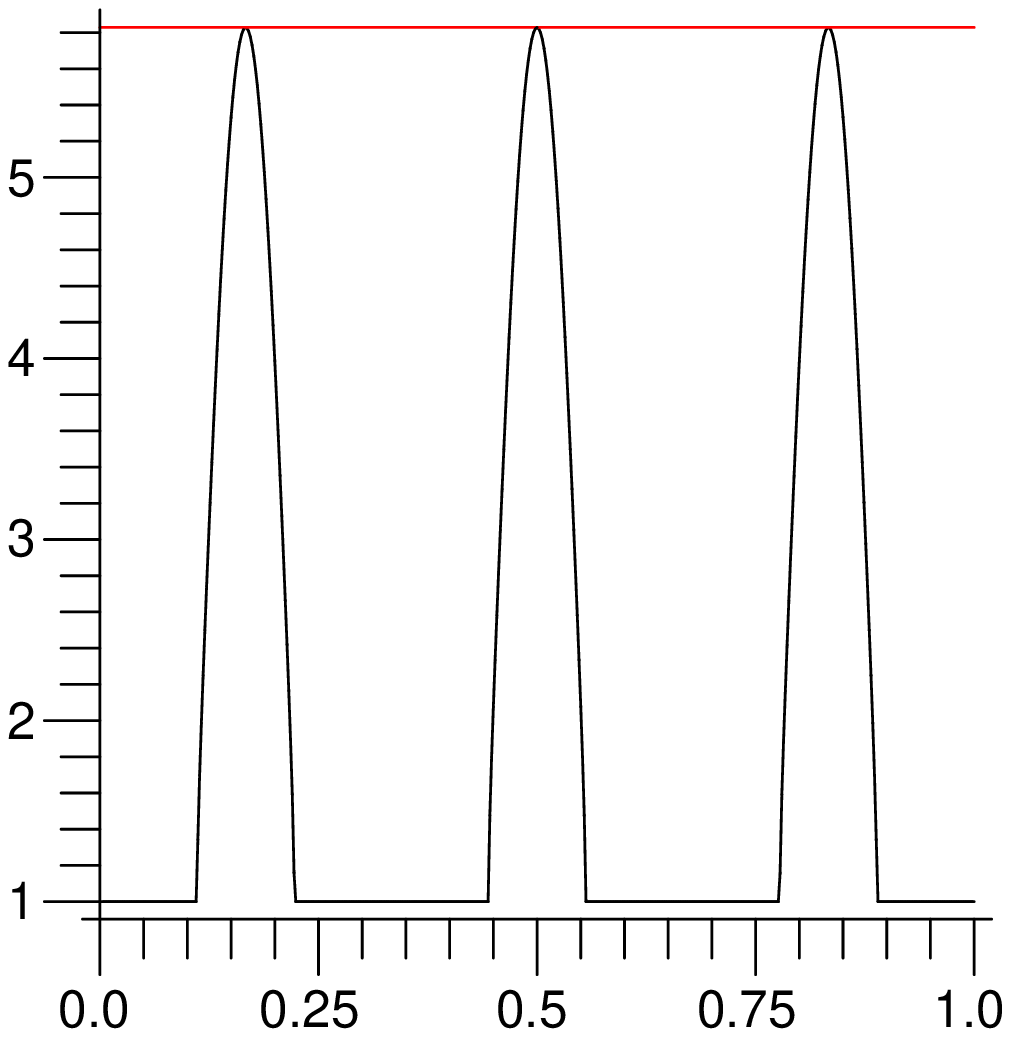}
      }
      &
      \scalebox{0.45}{
        \includegraphics[bb=80 250 520 575,clip]{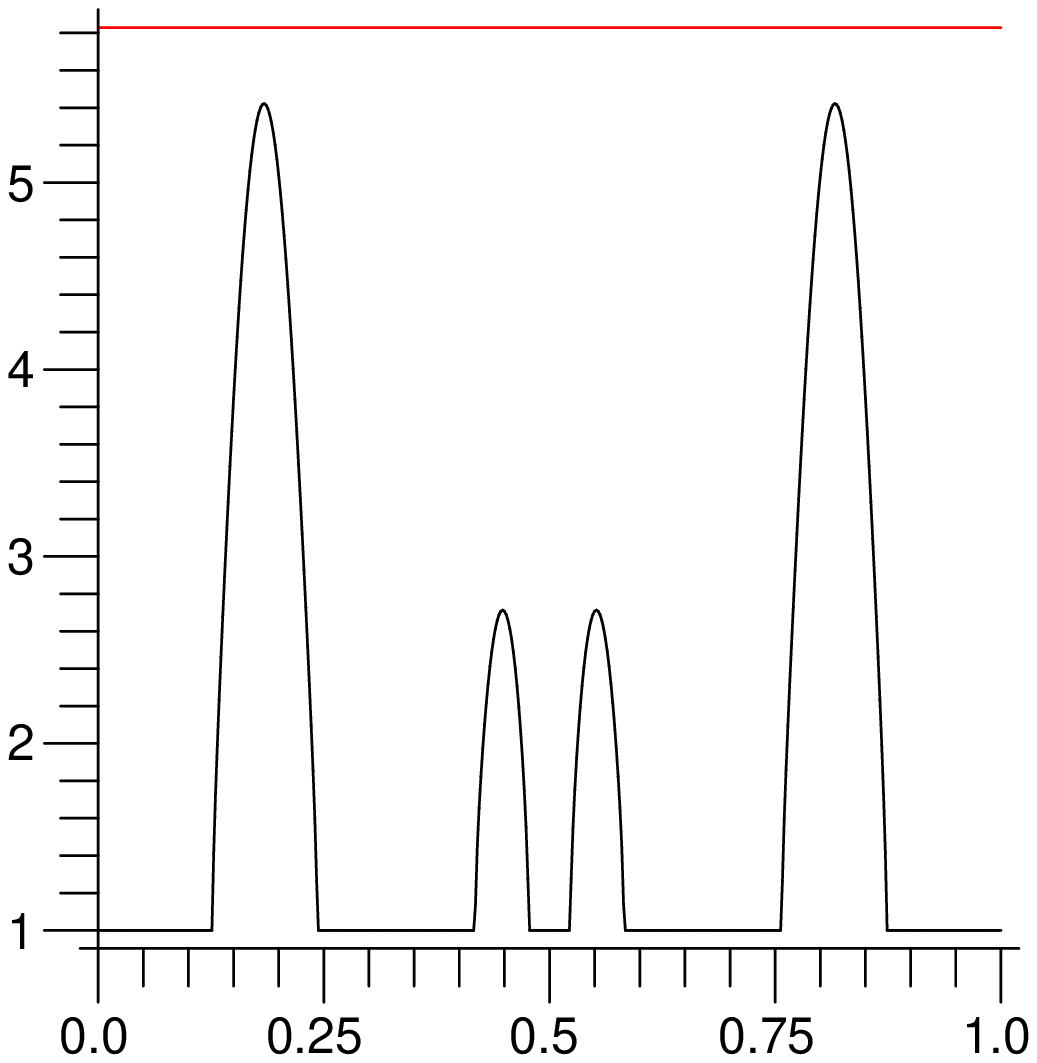}
      } \\
      (a) & (b)
    \end{array}
  \end{equation*}
  \caption{The Burau estimate (a) for~$\beta'_1$ and (b)
    for~$\beta'_2$.  The horizontal lines represent the growth
    rate~$\lambda \sim 5.828$.}
  \label{fig:specrad2}
\end{figure}
\end{example}

\begin{example}
  Let~$\beta'' \in B_{10}$ be the braid
  \begin{equation*}
    \beta'' = \beta_3 \cdot \sigma_{1,9,1} \cdot \sigma_{1,1,9}
  \end{equation*}
  where~$\beta_3 \in B_9 \subset B_{10}$ is as in Example \ref{ex:1}
  above\footnote{As is usual we regard~$B_9$ as the subgroup
    of~$B_{10}$ generated by the first eight generators~$\sigma_1,
    \dots, \sigma_8$.}.  See Figure \ref{fig:areduciblebraid2}(a).
  The Thurston normal form~$\red$ of~$\beta''$ has one pseudo-Anosov
  component~$\red_1 = \red|_{X_1}$, topologically conjugate to the
  pseudo-Anosov component of the Thurston normal form of~$B_3$.  By
  Theorem \ref{theorem:spectral radius for reducible braids} this
  component is Burau orientable and the Burau estimate is sharp at
  each of the cubic roots of~$-1$.  However,~$\beta''$ is
  distinguished from~$\beta_3$ by the fact that the eigenvalue
  of~$B(-e^{2 \pi i j/3})$ of modulus~$\lambda$ is not always real
  (see Figure \ref{fig:areduciblebraid2}(b)).  This can be seen as
  follows.  Lemma \ref{lemma: reducible vs orientability} below shows
  that for~$k$ a multiple of~$6$ the lift of~$\red_1$ to the~$k$-fold
  cover~$\tk{X}_1$ has orientable foliations.  Because~$\red$ acts as
  a full twist in the outer component of~$D$, the three connected
  components of~$\tk{X}_1$ are permuted cyclically by~$\tk{\red}$.  So
  by Theorem \ref{Main Algebra Theorem} and Lemma \ref{lemma:
    orientability vs eigenvalues}, the matrices~$B(\omega)$,
  for~$\omega$ a~$k$-th root of unity, contribute to~$\tk{\red}_*$
  exactly three eigenvalues of modulus~$\lambda$, differing from each
  other by the cubic roots of unity.  These eigenvalues are the
  extremal points of the curve in Figure
  \ref{fig:areduciblebraid2}(b).
  \begin{figure}
    \begin{center}
      \begin{equation*}
        \begin{array}{c@{\hspace{0.5in}}c}
        \put(80,-20){(a)}
        \put(310,-20){(b)}
          \scalebox{0.7}{
            \includegraphics{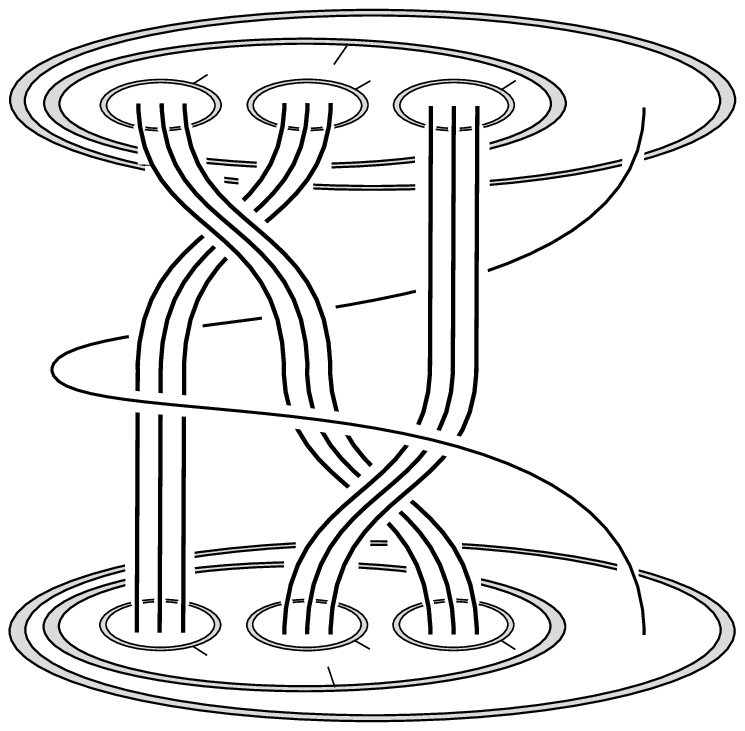}
          }
          &
          \scalebox{0.4}{
            \includegraphics[bb=0 230 575 640,clip]{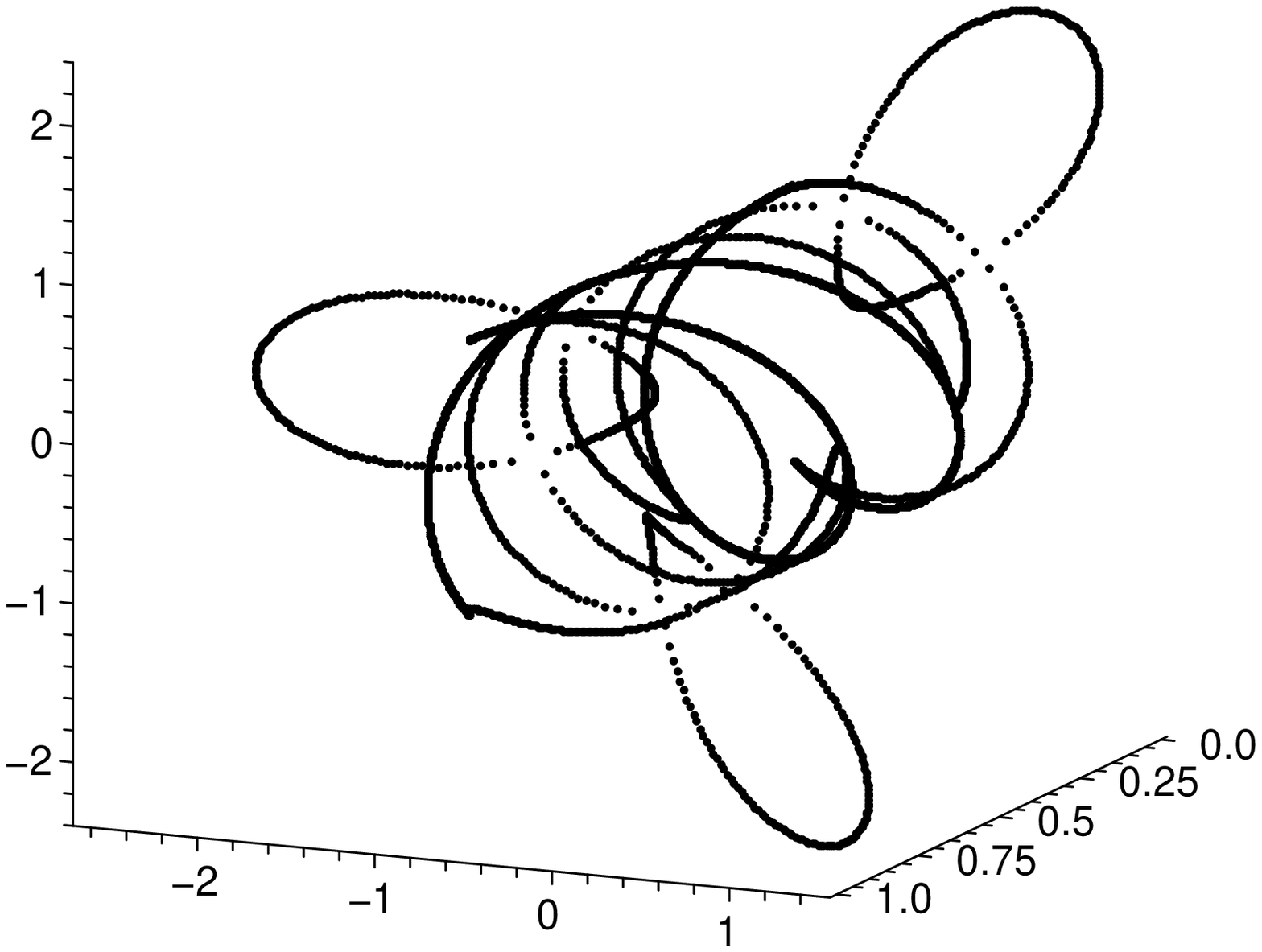}
          }
        \end{array}
      \end{equation*}
      \caption{(a) The braid~$\beta''$; (b) the eigenvalues of~$B(e^{2
          \pi i \theta})$ for~$\theta \in [0,1]$.  Here the rightmost
        axis represents the unit interval~$[0,1]$ and the other two
        axes represent the complex plane.}
    \label{fig:areduciblebraid2}
  \end{center}
  \end{figure}
\end{example}

\subsection{Proof of Theorem \ref{theorem:spectral radius for
  reducible braids}}
We will prove Theorem \ref{theorem:spectral radius for reducible
  braids} by a series of lemmas.  The main idea of the proof, already
suggested by the above examples, is the following.  As we lift a
pseudo-Anosov component~$\red_i$ of the Thurston normal form~$\red$ of
the chosen braid~$\beta$ into successive covering spaces~$\tk{D}$, the
supporting surface of the lift of~$\red_i$ may become disconnected.
If the lifted invariant foliations also become orientable, then Lemma
\ref{lemma: orientability vs eigenvalues} implies the existence of
several eigenvalues of modulus~$\lambda$ (counted with geometric
multiplicity).  By Theorem \ref{Main Algebra Theorem} these
eigenvalues are distributed among those of the matrices~$B(\omega)$,
where~$\omega$ is an~$k$-th root of unity.  An analysis of the action
of~$\tinf{\red}$ on the first homology of~$\tinf{X}_i$ then shows that
it is precisely the~$a_i$-th roots of~$-1$ which contribute such an
eigenvalue.

We begin with a lemma which provides a stronger version of the first
part of Theorem \ref{theorem:spectral radius for reducible braids}.
Let~$\red = \red_1 \cup \cdots \cup \red_m$ be the Thurston normal
form of the braid~$\beta \in B_n$, let~$\red_i$ be a pseudo-Anosov
component of~$\red$, and let~$\cF^u$ and~$\cF^s$ denote the invariant
foliations of~$\red_i$.  For~$k \geq 1$ we consider the~$k$-fold Burau
covering space~$\tk{p}:\tk{D}\raw D$ and write~$\tk{X}_i =
(\tk{p})^{-1}(X_i)$ and~$\tk{\red}_i = \tk{\red}|_{\tk{X}_i}$.  As in
the definition of~$a_i$ in \eqref{eq:a_i}, let~$X_{i0}$ be a connected
component of~$X_i$, and let~$x_1, \dots, x_{r'}$ be the punctures
in~$X_{i0}$ and~$O_{r'+1}, \dots, O_r$ the deleted disks.  Further,
for each~$1 \leq j \leq r$ let~$\kappa_j$ be the order of the
singularity which~$\cF^u$ and~$\cF^s$ exhibit at~$x_j$ (if~$j \leq
r'$) or on~$\bd O_j$ (if~$r' < j \leq r$).  We have:
\begin{lemma}
  \label{lemma: reducible vs orientability}
  Let~$\red_i = \red|_{X_i}$ be a pseudo-Anosov component of the
  Thurston normal form~$\red$ of~$\beta$, and let~$a_i$ be defined as
  in \eqref{eq:a_i}.  Write~$a = a_i$ and suppose that~$a = 2^u a'$
  where~$a'$ is odd.  Then the following are equivalent:
  \begin{tightlist}
  \item{(1)} $\red_i$ 
  is Burau orientable, i.e.\ there exists~$k \geq 1$ such that
    the lifts of~$\cF^u$ and~$\cF^s$ to~$X_i^{(k)}$ are orientable;
  \item{(2)} The set of~$k$ for which the lifts of~$\cF^u$ and~$\cF^s$ to~$X_i^{(k)}$ are orientable is the set
    of multiples of~$2^{u+1}$.
  \item{(3)} All singularities of~$\cF^u$ and~$\cF^s$ in the interior of~$X_{i0}$ have
    even order, and for all~$1 \leq j \leq r$,~$m_j / a \equiv \kappa_j \mod
    2$.  In other words~$m_j$ is an odd (even) multiple of~$a$
    whenever~$\kappa_j$ is odd (even).
  \end{tightlist}
\end{lemma}
\begin{proof}
  The implication~$(2) \implies (1)$ is trivial; we will prove~$(1)
  \implies (3)$ and~$(3) \implies (2)$.  Except for an application of
  Lemma \ref{orientlemma}, the proof is essentially algebraic -- it
  involves only the relevant morphisms of~$H_1(X_{i})$ into~$\Z$
  and~$\Z_k$.  For~$1 \leq j \leq r$ let~$\alpha_j$ be the homology
  class of a small clockwise loop around the puncture~$x_j$ (if~$j
  \leq r'$) or around the disk~$O_j$ (if~$j > r'$).
  Then~$H_1(X_{i0})$ can be identified with the subgroup of~$H_1(D)$
  generated by the~$\alpha_j$'s.  By definition~$\tau(\alpha_j) = m_j$
  for each~$j$, and~$\tau_k(\alpha_j) = \xi_k \tau(\alpha_j) = m_j
  \mod k$, where~$\xi_k:\Z\raw\Z_k$ is the quotient homomorphism.
  When all singularities of~$\cF^u$ and~$\cF^s$ in the interior
  of~$X_i$ have even order, we let~$\orient{\rho}:H_1(X_{i})\raw \Z_2$
  be the morphism associated to the orientation cover of~$\cF^u$
  and~$\cF^s$.  Because~$\red_i$ preserves the foliations and permutes
  the components cyclically,~$\orient{\rho}$ is determined by its
  values on~$H_1(X_{i0})$: namely~$\orient{\rho}(\alpha_j) = \kappa_j
  \mod 2$ for all~$j$.
  
  According to Lemma \ref{orientlemma} we know that for~$k \geq 1$ the
  following two statements are equivalent:
  \begin{enumerate}
  \item[(i)]~The lifts of~$\cF^u$ and~$\cF^s$ to~$X_i^{(k)}$ are orientable;
  \item[(ii)]~All singularities in the interior of~$X_{i}$ have even order,
    and there exists a homomorphism~$\delta_k:\im\tau_k|_{H_1(X_i)}\raw\Z_2$
    such that~$\delta_k\circ\tau_k=\orient{\rho}$ on~$H_1(X_i)$.
  \end{enumerate}
  In addition, since~$\tau_k \red_* = \tau_k$, the
  formula~$\delta_k\circ\tau_k=\orient{\rho}$ holds on~$H_1(X_i)$ whenever it
  holds on~$H_1(X_{i0})$.

  We note some elementary facts about cyclic groups.  
  For~$q \geq 1$ odd, there is no
  nontrivial homomorphism~$\Z_q \mapsto \Z_2$, while for~$q$ even there is
  only one such homomorphism: namely, the homomorphism which sends odd
  multiples of the generator to~$1$ and even multiples to~$0$ (a property
  which is independent of the choice of generator).  And if~$q$ is odd then
  every subgroup of~$\Z_q$ has odd order.

  Suppose that~$(1)$ holds, that is, that~$\red_i$ is pseudo-Anosov
  and for some~$k > 0$ the lifts of~$\cF^u$ and~$\cF^s$ to~$\tk{X}$
  are orientable.  Since~$\tk{p}:\tk{X}_i\raw X_i$ is ramified only
  around the punctures, this implies that all singularities in the
  interior of~$X_i$ have even order.  Since~$\im \tau|_{H_1(X_i)} = a
  \Z$, we have~$\im \tau_k|_{H_1(X_i)} = \xi_k(a\Z) = \xi_k(a) \Z_k$.
  The Euler-Poincare-Hopf formula \eqref{EPHFormula} shows
  that~$\cF^u$ and~$\cF^s$ have at least one~$1$-pronged singularity
  in~$X_{i0}$, which must lie at a puncture or on the boundary of a
  deleted disk.  Thus~$\orient{\rho}$ is onto, so from the previous
  paragraph we see that~$(ii)$ holds precisely when~$k$ is even and
  for all~$1 \leq j \leq r$ we have:~$\tau_k(\alpha_j)$ is an odd
  (even) multiple of~$\xi_k(a)$ iff~$\kappa_j$ is odd (even).  If~$k$
  is even then~$\xi_k$ sends odd (even) multiples of~$a$ to odd (even)
  multiplies of~$\xi_k(a)$ and so this implies the same statement
  in~$\Z$: that is, that~$(3)$ holds.  Thus~$(1) \implies (3)$.

  On the other hand suppose~$(3)$ holds and set~$k=2^{u+1}$.  Since~$a = 2^u
  a'$ where~$a'$ is odd, we see that~$\im \tau_k|_{H_1(X_i)}$ is the
  subgroup~$\{0,\xi_k(a)\}$ of~$\Z_k$ having order two.  Let~$\delta_k$ be the
  isomorphism of this subgroup onto~$\Z_2$.  Again, since~$k$ is even, to say
  that~$\xi_k(m_j) = \xi_k(a)$ is to say that~$m_j$ is an odd multiple of~$a$.
  By~$(3)$ this happens precisely when~$\kappa_j$ is odd; hence~$\delta_k
  \tau_k = \orient{\rho}$ on~$H_1(X_i)$.  Therefore since~$(i) \Leftrightarrow
  (ii)$ the lifts of~$\cF^u$ and~$\cF^s$ to~$X_i^{(k)}$ are orientable.
  Therefore the set of~$k$ for which the lifted foliations in~$X^{(k)}$ are
  orientable contains~$2^{u+1}$, and, it follows, also all multiples
  of~$2^{u+1}$.

  To complete the proof that~$(3) \implies (2)$, suppose for a contradiction
  that~$(3)$ holds and that for some~$k$ which is not a multiple of~$2^{u+1}$,
  the lifts of~$\cF^u$ and~$\cF^s$ to~$X_i^{(k)}$ are orientable.  The last
  condition implies that~$k$ is even, so we must have~$u > 0$.  Let~$\delta_k$
  be the homomorphism supplied by~$(ii)$.  Since~$\orient{\rho}$ is
  onto,~$\delta_k$ must send the generator~$\xi_k(a)$ of~$\im
  \tau_k|_{H_1(X_i)}$ to~$1 \in \Z_2$.  Now write~$k = 2^{u'} k'$ where~$k'$
  is odd and~$u' \leq u$.  The order of~$\xi_k(a)$ in~$\Z_k$ is given
  as~$\zeta = \lcm(a,k) / a = \lcm(a',k')/a'$, which is odd.  But then~$0 =
  \delta_k(0) = \delta_k (\zeta \xi_k(a)) = \zeta \delta_k \xi_k(a) = \delta_k
  \xi_k(a) = 1$, a contradiction.  Therefore no such~$k$ can exist.  This
  completes the proof.  \QED
\end{proof}

To prove the second part of Theorem \ref{theorem:spectral radius for
  reducible braids}, we will require two more lemmas.  The first of
these expresses~$H_1(\tinf{D})$ in terms of the first homology groups
of the~$\tinf{X}_i$'s.  Since each~$\tinf{X}_i$
is~$\tinf{\red}$-invariant, we will be able to use this lemma to
factorize the characteristic polynomial of the reduced Burau
matrix~$B(\beta)$.

\begin{lemma}
  \label{lemma:reducible homology splitting}
  Let~$\red = \red_1 \cup \dots \cup \red_m \in \beta$ be as above and
  let~$\tinf{\red}$ be the lift of~$\red$ to~$\tinf{D}$.  Let~$\tinf{X}_i =
  (\tinf{p})^{-1}(X_i)$ be the pull-back to~$\tinf{D}$ of the supporting
  surface~$X_i$ of~$\red_i$.  Then, as an Abelian group,~$H_1(\tinf{D})$
  decomposes as a direct sum of subgroups
  \begin{equation}
    \label{eq: reducible homology splitting}
    H_1(\tinf{D}) = \left(\bigoplus_{i=1}^m H_1(\tinf{X}_i)\right) \oplus V
  \end{equation}
  in which each~$H_1(\tinf{X}_i)$ is~$\tinf{\red}_*$-invariant, and~$V$ is a
  free abelian group of finite rank.
\end{lemma}

\begin{proof}
  We use the Mayer-Vietoris sequence inductively.  Let~$\tinf{A} =
  (\tinf{p})^{-1}(A)$ be the set of points of~$\tinf{D}$ covering the reducing
  annuli in~$D$.  For~$1 \leq i \leq m$ define~$A_i$ to be the set of those
  connected components of~$\tinf{A}$ whose boundary intersects~$\tinf{X}_i$, but
  does not intersect~$\tinf{X}_{j}$ for any~$1 \leq j < i$.  Then for each~$i$
  we let~$Z_i = \tinf{X}_i \cup A_i$.  Also, for~$0 \leq i \leq m$ let~$Y'_i =
  \bigcup_{j > i} Z_j$ and when~$i > 0$ let~$Y_i = Y'_i \cup A_i$.  These spaces
  satisfy
  \begin{align*}
    Z_i \cap Y_i = A_i,
    \quad & \qquad
    Z_i \cup Y_i = Y'_{i-1}  \\
    \intertext{for each~$1 \leq i \leq m$, and}
    Y'_{0} = \tinf{D},
    \quad & \qquad 
    A_m = Y'_m = Y_m = \emptyset.
  \end{align*}
  Since~$\tinf{X}_i$ is~$\tinf{\red}$-invariant for each~$i$, so too are~$A_i,
  Z_i, Y_i$ and~$Y'_i$.

  For any reducing annulus~$a \in A$, each boundary component of~$a$ is also the
  boundary of some (possibly punctured) subdisk~$D'$ of~$D$.  In fact, since no
  component of~$\bd a$ is allowed to be null-homotopic,~$D'$ must have at least
  one puncture and so we will have~$\tau([\bd D']) \neq 0$.  It follows that any
  connected component of~$\tinf{A}$ covering~$a$ is an infinite strip,
  isomorphic as a covering space to the universal cover of~$a$.  In
  particular~$H_1(A_i) = 0$ for each~$i$.

  Therefore, regarding~$Z_i$ and~$Y_i$ as subspaces of their union~$Y'_{i-1}$,
  the Mayer-Vietoris sequence for the pair~$(Z_i,Y_i)$ contains the segment
  \begin{equation}
    \label{eq: our mayer-vietroris sequence}
    \xymatrix{
      0 \ar[r]
      &
      {H_1(Z_i) \oplus H_1(Y_i)} \ar[r]
      &
      H_1(Y'_{i-1}) \ar[r]^{\bd}
      &
      H_0(A_i) \ar[r]
      & \dots
    }
  \end{equation}
  All groups in the above sequence are free abelian.  Noting that~$\tinf{X}_i$ is a
  deformation retract of~$Z_i$ and that~$Y'_{i-1}$ is a deformation retract
  of~$Y_{i-1}$, we have
  \begin{equation*}
    H_1(Y_{i-1}) \isomorphic H_1(\tinf{X}_i) \oplus H_1(Y_i) \oplus V_i
  \end{equation*}
  where~$V_i$ is canonically isomorphic to~$\im \bd \subset H_0(A_i)$.
  Since~$A_i$ has finitely many connected components,~$H_0(A_i)$ and~$V_i$ each
  have finite rank.  The lemma is proved by applying this formula inductively,
  starting with~$H_1(Y_{m-1}) = H_1(\tinf{X}_m)$, and setting~$V =
  \bigoplus_{i=1}^m V_i$. \QED\
\end{proof}

\begin{corollary}
  \label{corollary:burau matrix splitting}
  For each~$1 \leq i \leq m$ let~$\tinf{\red}_i$ be the restriction
  of~$\tinf{\red}$ to~$\tinf{X}_i$, and denote by~$g_i$ the action
  of~$\tinf{\red}_i$ on~$H_1(\tinf{X}_i)$.  As in Section \ref{sect:finite},
  regard~$H_1(\tinf{D})$ and~$H_1(\tinf{X}_i)$ as modules over the ring~$R =
  \Z[t^{\pm 1}]$.  Then~$\bigoplus_{i=1}^m H_1(\tinf{X}_i)$ is a submodule
  of~$H_1(\tinf{D})$ of full rank.  Moreover, the characteristic polynomial of
  the module isomorphism~$\tinf{h}_*:H_1(\tinf{D})\raw H_1(\tinf{D})$ is given
  by
  \begin{equation*}
    \chi(\tinf{h}_*) = \prod_{i=1}^m \chi(g_i).
  \end{equation*}
\end{corollary}
\begin{proof}
  Let~$G = \bigoplus_{i=1} H_1(\tinf{X}_i)$.  Because the~$R$-module structure
  on~$H_1(\tinf{D})$ and on~$H_1(\tinf{X}_i)$ are each defined using the same
  group of deck transformations acting respectively on~$\tinf{D}$ and
  on~$\tinf{X}_i$, they coincide, making~$G$ a submodule of~$H_1(\tinf{D})$.
  
  As an~$R$-module,~$H_1(\tinf{D})$ is free and of rank~$n-1$.  We claim that~$G
  \subseteq H_1(\tinf{D})$ also has rank~$n-1$.  For suppose~$\rank(G) < n-1$.  We
  may choose~$w \in H_1(\tinf{D})$ so that~$G \cup \{w\}$ spans a submodule
  of~$H_1(\tinf{D})$ of rank strictly larger than~$\rank(G)$.
  Since~$H_1(\tinf{D})$ is free, it follows that~$w + G$ generates a
  free~$R$-submodule of~$H_1(\tinf{D}) / G$.  But in the notation of Lemma
  \ref{lemma:reducible homology splitting}, we have~$H_1(\tinf{D}) / G
  \isomorphic V$ as Abelian groups.  Therefore~$H_1(\tinf{D}) / G$ has finite
  rank as an Abelian group, a contradiction.

  By definition~$\tinf{h}_*|_G = \bigoplus_{i=1}^m g_i$.  Because~$G$ is a
  submodule of~$H_1(\tinf{D})$, the characteristic polynomial of~$\tinf{h}_*|_G$
  divides the characteristic polynomial of~$\tinf{h}_*$ on~$H_1(\tinf{D})$.
  Since~$G$ has full rank both polynomials have the same leading
  coefficient~$x^{n-1}$ and so they are equal.  This completes the proof. \QED\
\end{proof}

The next lemma shows how the characteristic polynomial~$\chi(g_i)$
of~$g_i$ reflects the permutation induced by~$\tinf{\red}_i$ on the
set of connected components of~$\tinf{X}_i$.

\begin{lemma}
  \label{lemma: explicit form for char poly of lifted map}
  Let~$\red_i = \red|_{X_i}$ be a component of the Thurston normal
  form~$h$ of~$\beta$, and let~$\tinf{\red_i}$ be the lift of~$\red_i$
  to~$\tinf{X}_i$.  Let~$g_i$ denote the action of~$\tinf{\red_i}$
  on~$H_1(\tinf{X}_i)$.  Suppose~$X_i$ has~$\ell_i$ connected
  components, and let~$d \in \Z$ be chosen so
  that~$(\tinf{\red})^{\ell_i}(Y) = T^d(Y)$ for any connected
  component~$Y$ of~$\tinf{X}_i$.  Then there exists an integer~$e > 0$
  such that the characteristic polynomial of~$g_i$ is of the form
  \begin{equation}
    \label{eq: char poly of lifted component}
    \chi(g_i) = x^{\ell_i e} + t^{d} P_1(t) x^{\ell_i (e-1)} + t^{2d} P_2(t) x^{\ell_i (e-2)} + \dots + t^{ed} P_e(t)
  \end{equation}
  where each~$P_j \in \Z[t^{\pm a_i}]$.
\end{lemma}
\begin{proof}
  As before, we let~$X_{i0}$ be a connected component of~$X_i$
  and~$\tinf{X}_{i0}$ its pull-back to~$\tinf{D}$.  Let~$e$ be the dimension
  of~$H_1(\tinf{X}_{i0})$ as a module over~$\Z[t^{\pm 1}]$.  Let~$\{\upsilon_j :
  1 \leq j \leq e\}$ be a basis for~$H_1(\tinf{X}_{i0})$.  We choose this basis
  so that all of the~$\upsilon_j$'s are represented by loops in one and the same
  connected component,~$Y$ say, of~$\tinf{X}_{i0}$.  Pushing forward this basis
  under each of the first~$\ell_i$ iterates of~$g_i$, we obtain a basis
  for~$H_1(\tinf{X}_{{i}})$.  With respect to this basis, the matrix of~$g_i$
  has block form
  \begin{equation*}
    g_i = \left(
      \begin{matrix}
        0 & & &  \dots & \Omega \\
        \id & 0 & &  & \\
        & \id & 0 &  & \vdots \\
        & & \ddots & \ddots & \\
        & & & \id & 0
      \end{matrix}
    \right),
  \end{equation*}
  having~$\ell_i \times \ell_i$ blocks each of dimension~$e \times e$,
  where~$\Omega$ is some matrix in~$\GL(e,\Z[t^{\pm 1}])$.

  The connected components of~$\tinf{X}_{i0}$ are in one to one
  correspondence with the cosets of~$a_i\Z$ in~$\Z$.
  Now~$(\tinf{\red})^{\ell_i}$ sends the chosen connected component~$Y
  \subset \tinf{X}_{i0}$ to the possibly different connected
  component~$T^{d} Y$ of~$\tinf{X}_{i0}$, where~$d \in \Z$ is
  well-defined up to adding a multiple of~$a_i$.  It follows that the
  entries of~$\Omega$ lie in the coset~$t^d \Z[t^{\pm a_i}]$
  of~$\Z[t^{\pm a_i}]$ in~$\Z[t^{\pm 1}]$.

  Any product of~$j$ elements of~$\Omega$ therefore lies in the coset~$t^{jd}
  \Z[t^{\pm a_i}]$ of~$Z[t^{\pm a_i}]$.  The lemma is proved when we note that
  the characteristic polynomial of the matrix of~$g_i$ is of the form
  \begin{equation*}
    x^{\ell_i e} + \omega_1 x^{\ell_i (e-1)} + \dots + \omega_e
  \end{equation*}
  where for each~$1 \leq j \leq e$,~$\omega_j$ is a sum of products of~$j$
  elements of~$\Omega$, so that~$\omega_j \in t^{jd} \Z[t^{\pm a_i}]$.  \QED
\end{proof}

\begin{corollary}
  \label{corollary: roots of char poly}
  The roots of~$\chi(g_i)$ are naturally divided into sets each of~$\ell_i a_i$
  roots.  If~$(x,t)$ and~$(x',t')$ are two roots of~$\chi(g_i)$ in the same set,
  then~$(x',t') = (\mu x, \nu t)$, where~$\nu$ is some~$a_i$-th root of unity
  and~$\mu$ is some~$\ell_i$-th root of~$\nu^d$.
\end{corollary}
\begin{proof} Let~$(x,t) \in \C \times \C$ and let~$\mu$ and~$\nu$ be as in the
statement.  A quick calculation using \eqref{eq: char poly of lifted component}
shows that the value of~$\chi(g_i)$ at~$(\mu x, \nu t)$ is exactly~$\nu^{ed}$
times its value at~$(x,t)$.  In particular if~$(x,t)$ is a root so too is~$(\mu
x, \nu t)$.  \QED
\end{proof}

\noindent\textbf{Proof of Theorem~\ref{theorem:spectral radius for reducible braids}:}
Lemma \ref{lemma: reducible vs orientability} establishes the first
part of Theorem \ref{theorem:spectral radius for reducible braids}.
To prove the second part, recall that~$I$ is the set of~$1 \leq i \leq
m$ such that~$\red_i$ is pseudo-Anosov, Burau orientable, and
satisfies~$\htop(\red_i) = \htop(\red)$.  For~$k > 0$ we write~$\eta_k
= e^{2 \pi i/k}$.  Also, for each~$i$ we write~$\tk{\red}_i =
\tk{\red}|_{\tk{X}_i}$ and~$\tk{g}_i = \tk{\red}_*|_{H_1(\tk{X}_i)}$
as above.

Let~$i \in I$.  With notation as above, the number of connected components
of~$\tinf{X}_i$ is~$\ell_i a_i$.  Since~$\tau_k = \xi_k \tau$, the same
statement is true if we replace~$\tinf{X}_i$ by~$\tk{X_i}$, for any~$k$ which is
a multiple of~$a_i$.  

Suppose in fact that~$k$ is a multiple of~$2a_i$. By Lemma \ref{lemma:
  reducible vs orientability}, the invariant foliations
of~$\tk{\red}_i$ in~$\tk{X_i}$ are orientable.  The~$\ell_i a_i$
connected components of~$\tk{X_i}$ are permuted by~$\tk{\red}_i$ in
cycles all having the same number of components.  Suppose the number
of such cycles is~$L$, so that each cycle contains~$\ell_i a_i / L$
components.  Lemma \ref{lemma: orientability vs eigenvalues} applies
to each such cycle, and we conclude that~$\tk{g}_i$ has
exactly~$\ell_i a_i / L$ eigenvalues of modulus~$\lambda$, each having
geometric multiplicity~$L$.  Indeed~$\epsilon \lambda$ is one such
eigenvalue, where~$\epsilon \in \{1,-1\}$ is chosen according to
whether~$\tk{\red}_i$ preserves or reverses the orientation of its
unstable foliation.  Counted with geometric multiplicity, there are
exactly~$\ell_i a_i$ eigenvalues of modulus~$\lambda$.

Theorem \ref{Main Algebra Theorem} now implies that these~$\ell_i a_i$
eigenvalues are distributed among the eigenvalues of the matrices~$M(\eta_k^j)$,
where~$M$ denotes the matrix of~$g_i:H_1(\tinf{X}_i)\raw H_1(\tinf{X}_i)$ as in
Section \ref{sect:finite}.  Because~$M(\nu_k^j)$ is obtained by
substituting~$\nu_k^j$ into~$M$, an eigenvalue~$x$ of~$M(\eta_k^j)$ corresponds
to a root of~$\chi(g_i)$ of the form~$(x,\nu_k^j)$.  In particular,
setting~$k=2a_i$, we conclude that there is some~$0 \leq j_0 < 2a_i$ such
that~$(\epsilon \lambda, \eta_{2a_i}^{j_0})$ is a root of~$\chi(g_i)$.

Meanwhile, since by Lemma \ref{lemma: reducible vs orientability} the
lifts of~$\cF^u$ and~$\cF^s$ to~$X^{(a_i)}_i$ are not orientable, a
similar argument shows that no~$a_i$-th root of unity can occur in
such a root of~$\chi(g_i)$.  Therefore~$j_0$ must be odd.  By
Corollary \ref{corollary: roots of char poly} we now see that every
element of the set
\begin{equation*}
  J_i = \{(\mu \epsilon \lambda, \eta_{2a_i}^{j_0 + j}) \st \text{$j$ is even and and~$\mu$
    is an~$\ell_i$-th root of~$\eta_{2a_i}^{jd}$}\}
\end{equation*}
is a root of~$\chi(g_i)$, where~$d$ is as in Lemma \ref{lemma: explicit form for
  char poly of lifted map}.  Since~$j_0$ is odd, the roots of unity occuring in
elements of~$J_i$ are precisely the~$a_i$-th roots of~$-1$.

For~$k$ a multiple of~$2a_i$, again by Theorem \ref{Main Algebra Theorem} we see
that~$J_i$ accounts for all of the~$\ell_i a_i$ eigenvalues of~$\tk{g}_i$ of
modulus~$\lambda$.  Furthermore, every root of unity can be written as a~$k$-th
root of unity for some~$k$ which is a multiple of~$2a_i$.  It follows that~$J_i$
is precisely the set of roots~$(x,\nu)$ of~$\chi(g_i)$ such that~$\nu$ is a root
of unity and~$|x| = \lambda$.   

By Corollary \ref{corollary:burau matrix splitting}, the characteristic
polynomial of~$B(\beta)$ is the product of those of the~$g_i$'s.  Therefore the
set of roots of~$\chi(B(\beta))$ is the union of those of the~$g_i$'s.  For~$i
\in I$, we have just accounted for all of the roots~$(x,\nu)$ of~$\chi(g_i)$
with~$|x|=\lambda$ and~$\nu$ a root of unity.  For~$i \not \in I$ a similar
argument shows that~$\chi(g_i)$ can have no such roots.  This completes the
proof. \QED

\noindent\textbf{Proof of Theorem~\ref{introthm2b}:}
  Fix~$\beta \in B_n$ and let~$\red = \red_1
  \cup \dots \cup \red_m$ be its Thurston normal form.  If the Burau
  estimate \eqref{burest} is sharp for the substitution of some root
  of unity~$\eta$ into the Burau matrix of~$\beta$, then by Theorem
  \ref{theorem:spectral radius for reducible braids} the set~$I$ of
  indices for which~$\red_i$ is pseudo-Anosov, Burau orientable, and
  satisfies~$\htop(\red_i) = \htop(\red)$ must be nonempty.
  Let~$K(\beta)$ be the smallest value of~$k$ for which there
  exists~$0 < j < k$ with~$\lambda = \spec(B(e^{2\pi i j/k}))$.
  According to Lemma \ref{lemma: reducible vs
    orientability},~$K(\beta)$ is the infimum, over~$i \in I$, of
  twice the power of~$2$ occurring in the prime factorization of~$a_i$.

  Now fix~$i \in I$ and let~$X_i$ be the supporting surface
  of~$\red_i$ and~$X_{i0} \subset X_i$ a connected component.  Adopt
  the notation of Lemma \ref{lem:a_i} and the definition of~$a_i$.  In
  particular~$r$ denotes the number of punctures and disks deleted
  from~$X_{i0}$.  Because~$\red_i$ is pseudo-Anosov, we must have~$r
  \geq 3$ by the Euler-Poincar\'e-Hopf formula.  Since~$a_i \leq
  \inf_{j} m_j$ and~$\sum_j m_j \leq n / \ell_i$ by definition, we
  therefore have~$a_i \leq n / r \ell_i \leq n / 3$.  Therefore the
  power of~$2$ in~$a_i$ is also at most~$\tfrac{1}{3}n$.  Since this holds for
  all~$i \in I$, we get~$K(\beta) \leq \tfrac{2}{3}n$, proving Theorem
  \ref{introthm2b}.\QED

\begin{remark}
  We remark that if~$n = 3n'$ where~$n'$ is a power of two, then the
  braid~$\beta_{n'} \in B_n$ constructed after Theorem
  \ref{theorem:spectral radius for reducible braids}
  satisfies~$K(\beta_{n'}) = \tfrac{2}{3} n$.  Thus for some braids,
  the bound in Theorem \ref{introthm2b} is attained.
\end{remark}

\medskip\noindent\textbf{Acknowledgements:} We would like to thank
Toby Hall for his plotting software as well as for many helpful
conversations, and Jean-Luc Thiffeault for sharing his stimulating
unpublished work on the Burau representation and for pointing out that
the results in \cite{Song2002} overlap with those in our
Theorem~\ref{propn: pA burau estimate}.

\bibliography{bibliography}
\bibliographystyle{alpha}

\medskip

\begin{center}
\begin{tabular}{p{6cm} p{6cm} }
Gavin Band & Philip Boyland \\
Dept. of Mathematics & Dept. of Mathematics\\
   University of Liverpool & University of Florida\\
   Peach St. & Little Hall\\
   Liverpool  L69 7ZL & Gainesville, FL 32605-8105 \\
   U.K. & U.S.A.\\
   \verb!g.band@liv.ac.uk! & \verb!boyland@math.ufl.edu! \\
\end{tabular}
\end{center}
\end{document}